\documentclass{amsart}

\usepackage{graphicx}
\usepackage{here}

 \newtheorem{theorem}{Theorem}[section]
 \newtheorem{corollary}[theorem]{Corollary}

 \theoremstyle{remark}
 \newtheorem{remark}{Remark}[section]
 \newtheorem{example}{Example}[section]

\begin{document}

\markboth{Kazuhiro Ichihara and Gakuto Kato}{Two-bridge links and stable maps into the plane}

\title{Two-bridge links and stable maps into the plane}

\author{Kazuhiro Ichihara}
\address{Department of Mathematics, College of Humanities and Sciences, Nihon University, 3-25-40 Sakurajosui, Setagaya-ku, Tokyo 156-8550, JAPAN\\
ichihara.kazuhiro@nihon-u.ac.jp}

\author{Gakuto Kato}
\address{Graduate School of Integrated Basic Sciences, Nihon University, 3-25-40 Sakurajosui, Setagaya-ku, Tokyo 156-8550, JAPAN}

\maketitle 

\begin{abstract}
We give a visual construction of stable maps from the $3$-sphere into the real plane enjoying the following properties: 
the set of definite fold points coincides with a given two-bridge link and the map only admits restricted types of fibers containing two indefinite fold points. 
As a corollary, we determine the stable map complexities defined by Ishikawa and Koda for sufficiently complicated two-bridge link exteriors. 
\end{abstract}




\section{Introduction} 
For smooth manifolds $M$ and $N$, let $C^{\infty}(M,N)$ be the set of smooth maps from $M$ to $N$ with the Whitney topology.
A smooth map $f : M \to N$ is called a \textit{stable map} if there exists a neighborhood $U_{f}$ of $f$ in $C^{\infty}(M,N)$ such that, for any map $g$ in $U_{f}$, there are diffeomorphisms $\Phi : M \to M$ and $\phi: N \to N$ satisfying $g = \phi \circ f \circ \Phi^{-1}$.

For a stable map from a smooth $3$-manifold to the real plane $\mathbb{R}^{2}$, it is known that the local singularities are classified into three types: definite fold points, indefinite fold points and cusp points (for example, see \cite{Levin'65}).
Moreover, for a closed smooth 3-manifold, Levine showed in \cite{Levin'65} that the cusp points of any such stable map can be eliminated by smooth homotopy. 
Thus, in the rest of this paper, we will assume that all the stable maps from $3$-manifolds into $\mathbb{R}^{2}$ have no cusp points.

Let $M$ be a smooth 3-manifold and $f: M \to \mathbb{R}^{2}$ a stable map from $M$ to $\mathbb{R}^2$ without cusp points. 
Then the set $S_{1} (f)$ of the indefinite fold points of $f$ is a smooth l-dimensional sub-manifold of $M$ and $f |_{S_{1}(f)}$ is an immersion with normal crossings (transverse double points). 
See \cite{Saeki'96} for example. 
In \cite{Saeki'96}, Saeki proved that there exists a stable map $f : M \to \mathbb{R}^{2}$ such that $f (S_{1} (f))$ has no such double points if and only if $M$ is a graph manifold (see \cite{Saeki'96} for the definition of graph manifolds).
Due to this result, it is natural to consider a stable map $f$ such that $f (S_{1} (f))$ has non-empty crossings.
The \textit{singular fibers}, i.e. the components of the preimage of such double points by $f$, are classified into the \textit{type $\mathrm{I\hspace{-1.2pt}I^{2}}$} and the \textit{type $\mathrm{I\hspace{-1.2pt}I^{3}}$} (\cite{Levin'65}). 
We denote the set of type $\mathrm{I\hspace{-1.2pt}I^{2}}$ singular fibers of $f$ (resp. type $\mathrm{I\hspace{-1.2pt}I^{3}}$ singular fibers of $f$) by $\mathrm{I\hspace{-1.2pt}I^{2}} (f)$ (resp. $\mathrm{I\hspace{-1.2pt}I^{3}} (f)$). 

\begin{figure}[H]
{\unitlength=1cm
\begin{picture}(12.5,1.2)(0,0)
\put(4,-0.5){\includegraphics[height=3cm,clip]{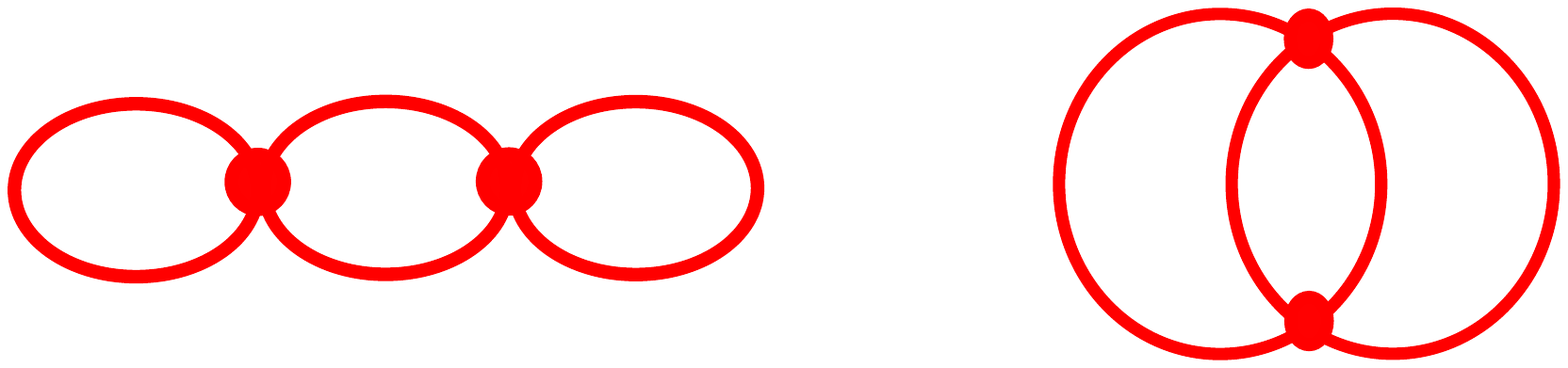}}
\put(4.5,0.2){type $\mathrm{I\hspace{-1.2pt}I^{2}}$}
\put(7,0.2){type $\mathrm{I\hspace{-1.2pt}I^{3}}$}
\end{picture}}
\caption{The types of the singular fibers.}
\label{sing.}
\end{figure}

In \cite{Ishikawa-Koda}, Ishikawa and Koda characterized the link $L$ in the $3$-sphere $S^{3}$ such that there exists a stable map $f : S^{3} \to \mathbb{R}^{2}$ satisfying that $L \subset S_{0} (f)$ and $(|\mathrm{I\hspace{-1.2pt}I^{2}} (f)|, |\mathrm{I\hspace{-1.2pt}I^{3}} (f)|) =(1,0)$. 
Here, $S_{0} (f)$ denotes the set of the definite fold points of $f$ and $|\mathrm{I\hspace{-1.2pt}I^{2}} (f)|$ and $|\mathrm{I\hspace{-1.2pt}I^{3}} (f)|$ denote the cardinalities of the sets $\mathrm{I\hspace{-1.2pt}I^{2}} (f)$ and $\mathrm{I\hspace{-1.2pt}I^{3}} (f)$, respectively.
Next, Furutani and Koda \cite{Furutani-Koda} characterized the link $L$ such that there exists a stable map $f : S^{3} \to \mathbb{R}^{2}$ satisfying that $L \subset S_{0} (f)$ and $(|\mathrm{I\hspace{-1.2pt}I^{2}} (f)|, |\mathrm{I\hspace{-1.2pt}I^{3}} (f)|) =(0,1)$.

In this paper, we show that, for some two-bridge link $L$, there exists a stable map $f : S^{3} \to \mathbb{R}^{2}$ such that $L=S_{0} (f)$ and $\mathrm{I\hspace{-1.2pt}I^{2}} (f)= \emptyset$ or $\mathrm{I\hspace{-1.2pt}I^{3}} (f)= \emptyset$.
Here, a non-trivial link in $S^3$ is called a \textit{two-bridge link} if it admits a diagram with exactly two maxima and two minima.

\begin{figure}[htb]
                \setlength\unitlength{1truecm}
                \begin{picture}(12.5,9.5)(0,0)
                    \put(-0.5,2){\includegraphics[height=9.5cm,clip]{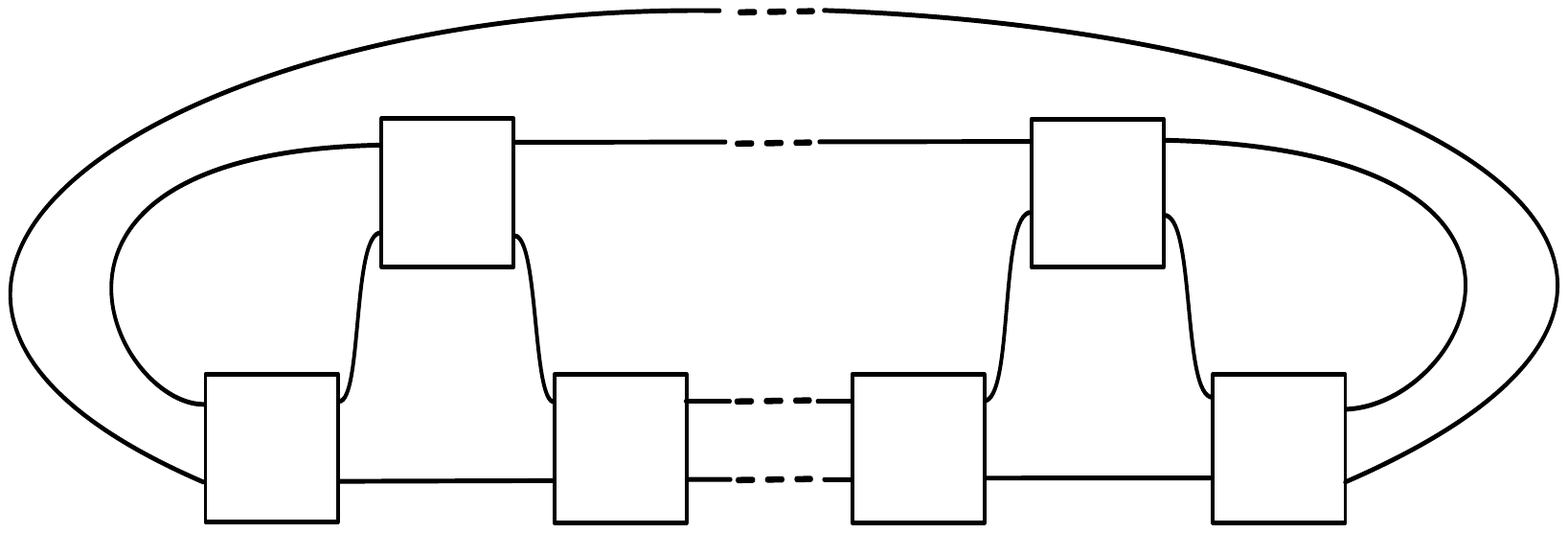}}
                    \put(2,5.7){$a_{1}$}
                    \put(3.4,7.8){$b_{1}$}
                    \put(4.8,5.7){$a_{2}$}
                    \put(7.1,5.7){$a_{m}$}
                    \put(8.5,7.8){$b_{m}$}
                    \put(9.8,5.7){$a_{m+1}$}

                    \put(0,-1){\includegraphics[height=6cm,clip]{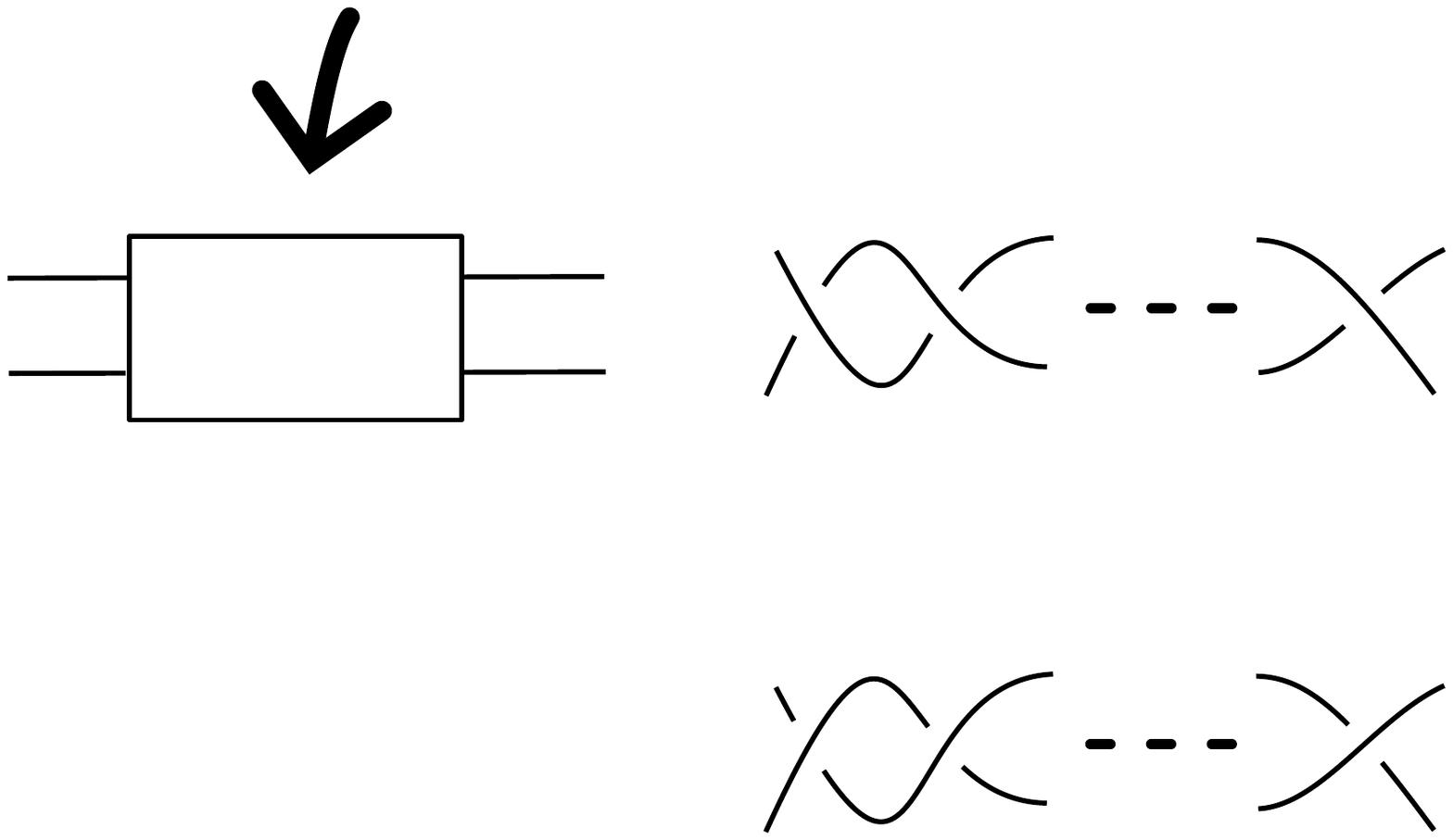}}
                    \put(4,3){$\rightsquigarrow$}
                    \put(1.3,3){$a_{i}$ or $b_{j}$}
                   
                    \put(8.5,3){if $a_{i}>0$ and $b_{j}<0$}
                    \put(8.5,0.7){if $a_{i}<0$ and $b_{j}>0$}
                \end{picture}
                \caption{Conway form $C(a_{1},b_{1},\cdots,a_{m},b_{m},a_{m+1})$. }
                \label{$36$}
            \end{figure}

For any two-bridge link $L$, there always exists a diagram of $L$ called the \textit{Conway form} $C(a_{1},b_{1},\cdots,a_{m},b_{m},a_{m+1})$ for $a_{i}, b_{j}\ne 0$ illustrated in Fig.~\ref{$36$}.
Then, our results are the following. 

    \begin{theorem}
        \label{main1}
            Let $L$ be a two-bridge link in $S^{3}$ represented by a Conway form $C(a_{1},b_{1},\cdots,a_{m},b_{m},a_{m+1} )$ with non-zero integers $a_{1},b_{1},\cdots,a_{m},b_{m},a_{m+1}$.
            If $b_{i}$ is even for all $i$, then there exists a stable map $f_{2}$ from $S^{3}$ into $\mathbb{R}^{2}$ satisfying that $S_{0}(f_{2})= L$, $|\mathrm{I\hspace{-1.2pt}I^{2}}(f_{2})| = 2m$ and $\mathrm{I\hspace{-1.2pt}I^{3}}(f_{2}) = \emptyset$.
    \end{theorem}

    \begin{theorem}
        \label{main2}
        Let $L$ be a two-bridge link in $S^{3}$ represented by a Conway form $C(a_{1},b_{1},\cdots,a_{m},b_{m},a_{m+1} )$ with non-zero integers $a_{1},b_{1},\cdots,a_{m},b_{m},a_{m+1}$. 
        If $b_{i}$ is even for all $i$, then there exists a stable map $f_{3}$ from $S^{3}$ into $\mathbb{R}^{2}$ satisfying that $S_{0}(f_{3})= L$, $\mathrm{I\hspace{-1.2pt}I^{2}}(f_{3}) = \emptyset$ and $|\mathrm{I\hspace{-1.2pt}I^{3}}(f_{3})| = \frac{1}{2} \sum_{i=1}^{m} | b_{i} |$.
    \end{theorem}

We remark that any two-bridge link of two components admits a Conway form $C(a_{1},b_{1},\cdots,a_{m},b_{m},a_{m+1} )$ with all the $b_{i}$'s are even. 
See \cite{Kawauchi-K.theory} for example. 
On the other hand, it is uncertain whether any two-bridge knot admits a Conway form $C(a_{1},b_{1},\cdots,a_{m},b_{m},a_{m+1} )$ of odd length with all the $b_{i}$'s are even.

Also, in \cite{Ishikawa-Koda}, Ishikawa and Koda introduced a complexity of a smooth 3-manifold, called the \textit{stable map complexity}, by using stable maps into $\mathbb{R}^2$. 
Precisely, the stable map complexity $\mathrm{smc}(M)$ of a compact 3-manifold $M$ with (possibly empty) boundary consisting of tori is defined as the minimal number of the weighted sums
$ \left| \mathrm{I\hspace{-1.2pt}I^{2}}(f) \right| + 2 \left| \mathrm{I\hspace{-1.2pt}I^{3}}(f) \right| $
for a stable map $f : M \to \mathbb{R}^2$. 

In \cite[Section 6]{Ishikawa-Koda}, the relationship between stable map complexities and hyperbolic volumes for hyperbolic 3-manifolds is discussed. 
Together with this, the following is obtained from Theorem~\ref{main1}. 

\begin{corollary}\label{cor}
Let $L$ be a two-bridge link which has the Conway form \\ $C(a_{1},b_{1},\cdots,a_{m},b_{m},a_{m+1})$ with positive even $a_{i}$, $b_{j}$ and $E(L)$ the exterior of $L$ which is obtained from $S^{3}$ by removing open tubular neighborhood of $L$. 
If $a_{i}$ and $b_{j}$ are sufficiently large, then $\mathrm{smc}(E(L)) = 2m$ holds. 
\end{corollary}


\section{A construction of $f_{2}$}

In this section, we give a construction of 
a stable map $f_{2}$ from $S^{3}$ into $\mathbb{R}^{2}$ satisfying that $S_{0}(f_{2})= L$, $|\mathrm{I\hspace{-1.2pt}I^{2}}(f_{2})| = 2m$ and $\mathrm{I\hspace{-1.2pt}I^{3}}(f_{2}) = \emptyset$ 
for a two-bridge link $L$ in $S^3$ when $L$ has a Conway form $C(a_{1},b_{1},\cdots,a_{m},b_{m},a_{m+1} )$ with all the $b_{i}$'s are even. 

\begin{figure}[H]
                \setlength\unitlength{1truecm}
                \begin{picture}(12.5,4)(0,0.7)
                    \put(0.2,-2){\includegraphics[width=1\textwidth,clip]{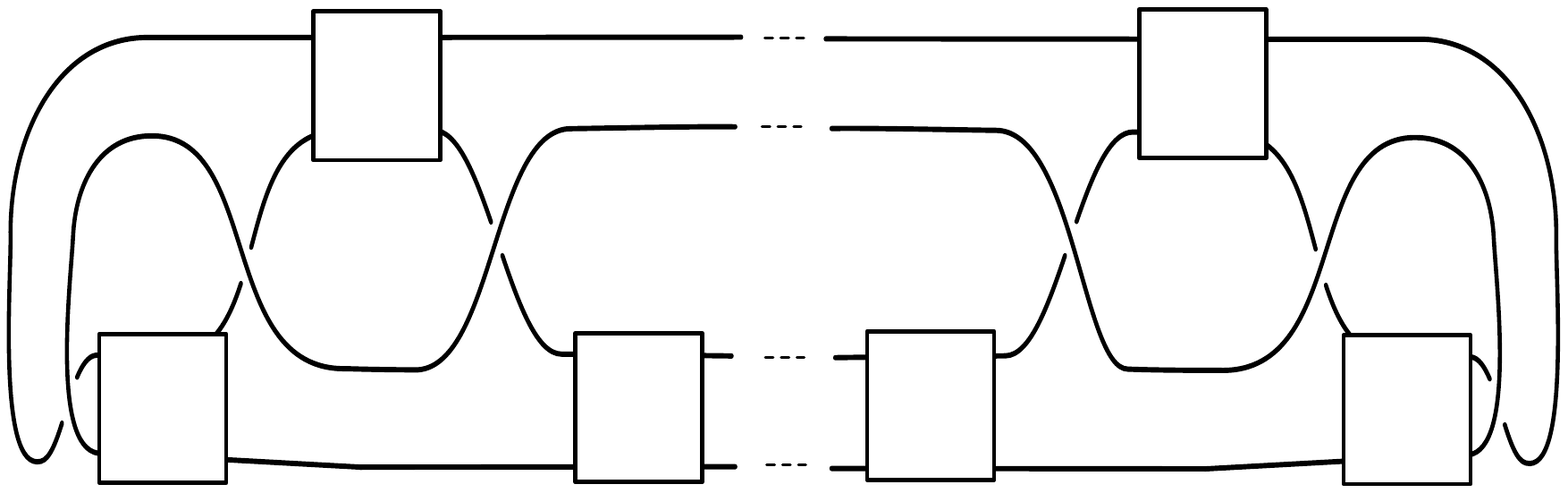}}

                    \put(1.4,1.2){$a_{1}$}
                    \put(5.25,1.2){$a_{2}$}
                    \put(7.5,1.2){$a_m$}
                    \put(11.2,1.2){$a_{m+1}$}
                    \put(3.1,3.7){$b_{1}$}
                    \put(9.7,3.7){$b_{m}$}
                \end{picture}
                \caption{The diagram $D^{\prime}$ of $L$.}
                \label{82}
            \end{figure}

First, by modifying the Conway form $C(a_{1},b_{1},\cdots,a_{m},b_{m},a_{m+1} )$ in Fig.~\ref{$36$}, we take another diagram $D^{\prime}$ of $L$ shown in Fig.~\ref{82}. 
            
            

\begin{figure}[H]
                \setlength\unitlength{1truecm}
                \begin{picture}(12.5,5)(0,0)
                    \put(1,-1){\includegraphics[width=0.8\textwidth,clip]{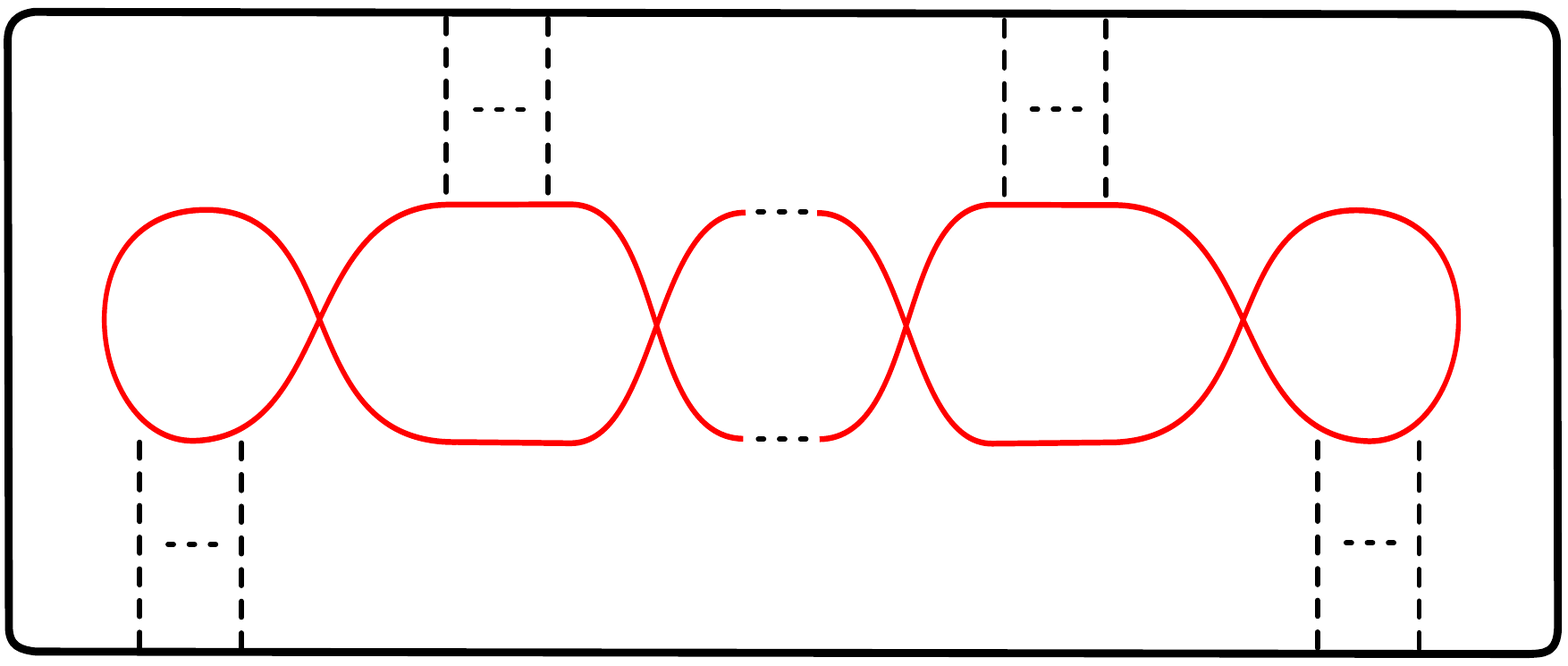}}
                    \put(3.75,4.7){\Large $\overbrace{ }^{|b_{1}|}$}
                    \put(7.35,4.7){\Large $\overbrace{ }^{|b_{m}|}$}
                    \put(1.75,0.35){\Large $\underbrace{ }_{|a_{1}|+1}$}
                    \put(9.1,0.35){\Large $\underbrace{ }_{|a_{m+1}|+1}$}
                \end{picture}
                \caption{The immersed curve $C^{\prime}$.}
                \label{90}
            \end{figure}

Next, we set a rectangular region $E$ on the plane $\mathbb{R}^{2}$ which contains $D^{\prime}$ inside. 
Let $C^{\prime}$ be an immersed curve obtained from $D^{\prime}$ as follows: 
Perform horizontal smoothing at all the double points adjacent to the outer region, removing the outermost circle appearing after the smoothing, and forgetting the crossing information of the remaining diagram. 
See Fig.~\ref{90}. 
In the figure, the immersed curve in the center depicts $C^{\prime}$ and the rectangular curve surrounding $C^{\prime}$ indicates the boundary of $E$. 
Also, vertical dotted lines are added to mark the crossings in the twists on $D^{\prime}$ associated to the $a_i$'s and the $b_{j}$'s. 
The horizontal dots between two vertical dotted lines indicate a set of vertical dotted lines. 

\begin{figure}[H]
                \setlength\unitlength{1truecm}
                \begin{picture}(12.5,6.6)(0,0.5)
                    \put(2,0.5){\includegraphics[width=0.75\textwidth,clip]{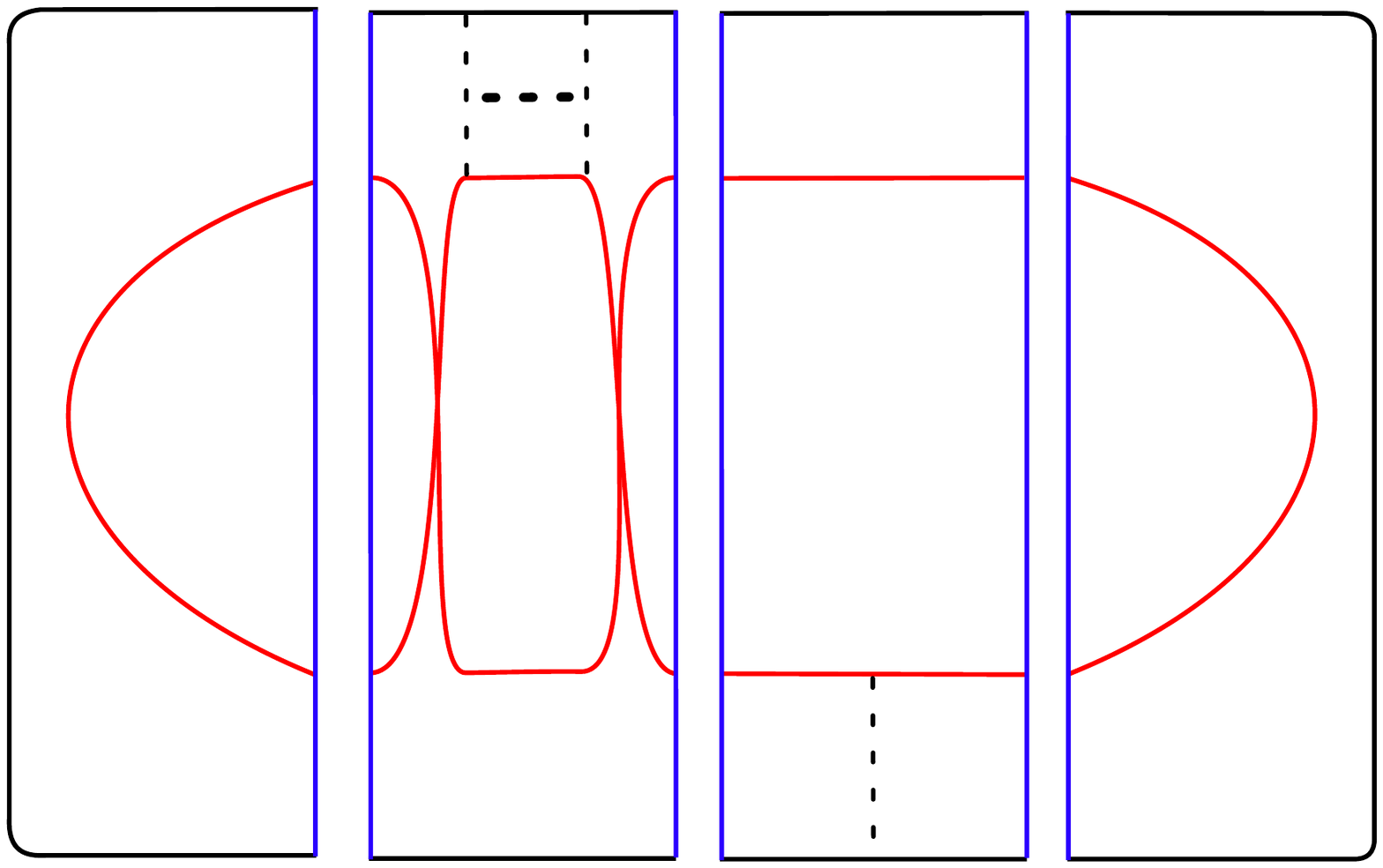}}
                    \put(5.15,6.8){\Large $\overbrace{ }^{|b_{i}|}$}
                    \put(2.6,0.5){Type $1$}
                    \put(5.1,0.5){Type $2$}
                    \put(7.5,0.5){Type $3$}
                    \put(9.8,0.5){Type $4$}
                \end{picture}
                \caption{The four types of rectangular regions.}
                \label{96}
            \end{figure}

Then, we take parallel segments $\gamma_{1},\gamma_{2},\cdots,\gamma_{n}$ on $E$ such that they separate $E$ into strips, say, $T_{0},T_{1},T_{2},\cdots,T_{n}$ so that $T_{k} \cap C^{\prime}$ looks like one of the four shapes depicted in Fig.~\ref{96} for each $k \in \{0,1,\cdots,n \}$. 
We refer to the rectangular regions separated by $\gamma_{k}$ as Type $1$, Type $2$, Type $3$ and Type $4$, as shown in Fig.~\ref{96}. 

\vspace{-0.7cm}

On the region of Type 2 in the figure, the horizontal dots between two vertical dotted lines indicate a set of vertical dotted lines. 
Note that the number of the regions of Type 2 is equal to $m$. 

\begin{figure}[H]
                \setlength\unitlength{1truecm}
                \begin{picture}(12.5,6.5)(0,1.2)
                    \put(0,0.5){\includegraphics[width=\textwidth,clip]{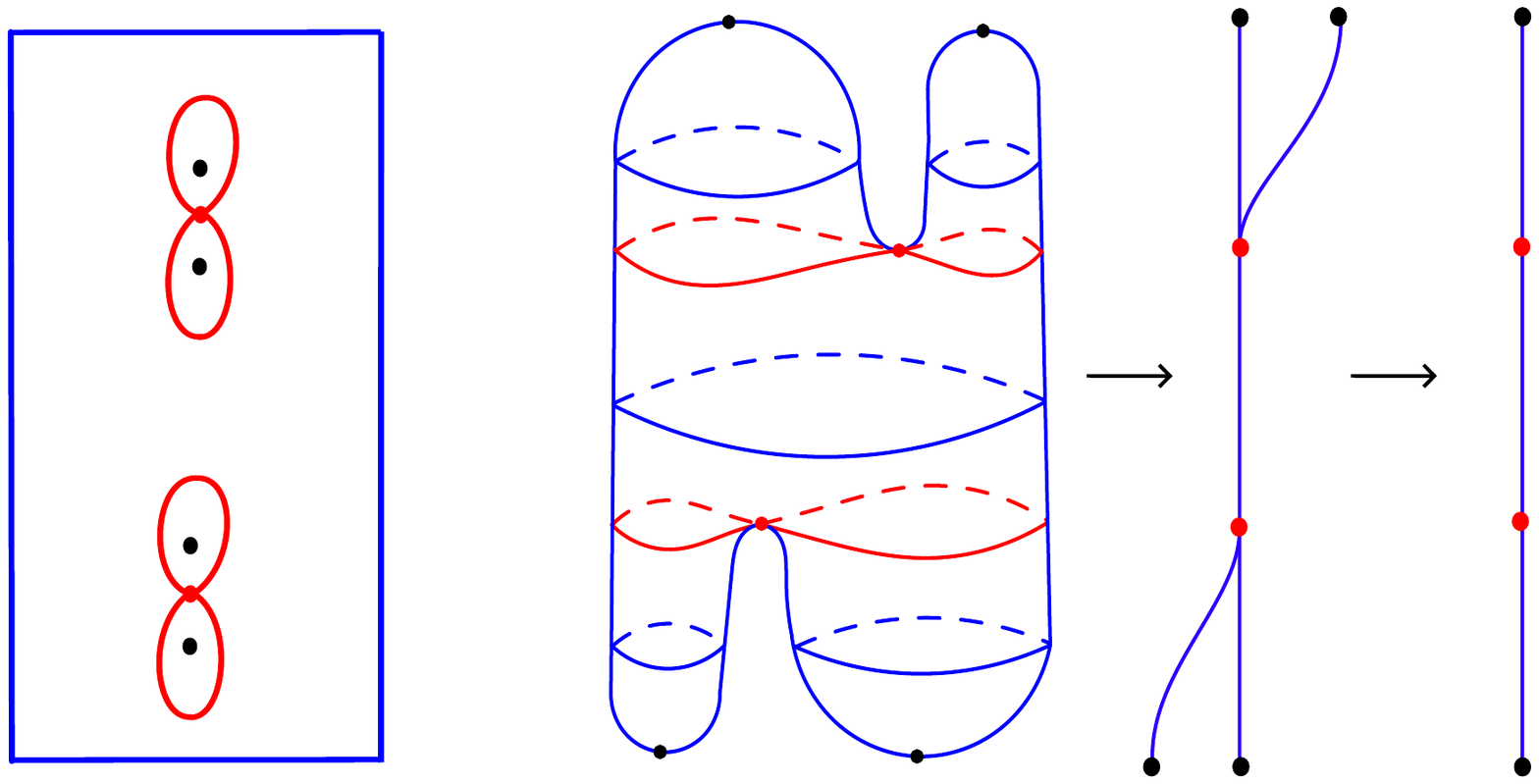}}
                    \put(3.6,4.9){\scalebox{3}{$\simeq$}}
                    
                    \put(1.7,6.8){$1$}
                    \put(1.6,5.5){$2$}
                    \put(1.6,3.7){$3$}
                    \put(1.5,2.3){$4$}

                    \put(5.8,8){$1$}
                    \put(8,8){$2$}
                    \put(5.4,1.5){$3$}
                    \put(7.4,1.5){$4$}

                    \put(10.1,8.1){$1$}
                    \put(10.9,8.1){$2$}
                    \put(9.3,1.4){$3$}
                    \put(10.2,1.4){$4$}

                    \put(12.2,8.1){$1,2$}
                    \put(12.2,1.4){$3,4$}
                \end{picture}
                \caption{A Morse function $\psi_{k} : F_{k} \to \gamma_{k} \subset E$.}
                \label{70-3}
            \end{figure}

Now, for each $k \in \{ 1,2,\cdots,n \}$, we construct a smooth map $\psi_{k} : F_{k} \to \gamma_{k}$ whose Reeb graph is the train track $\tau_{k}$ as shown in Fig.~\ref{70-3}, where $F_{k}$ is homeomorphic to the 2-sphere $S^{2}$. 
We regard $\psi_{k}$ as a Morse function on $F_{k} \simeq S^{2}$.

On the sphere $F_k$ in Fig.~\ref{70-3} (middle), the circles show regular fibers and the curves with singularity are singular fibers of $\psi_{k}$. 
The left figure exhibits a flattened image of $F_k$ on the plane. 

For a standard projection $\mathbb{R}^{3} \to \mathbb{R}^{2}$ to obtain $D^{\prime}$ from $L$, the preimage of a line on $\mathbb{R}^{2}$ containing $\gamma_{k}$ gives an embedded sphere in $S^{3}$. 
We identify it with $F_{k}$ for each $k  \in \{1,2,\cdots,n \}$, and decompose $S^{3}$ into $N_{0},N_{1},N_{2},\cdots,N_{n}$ by the spheres $F_{1},\cdots,F_{n}$. 
Here $N_{0}$ and $N_{n}$ are homeomorphic to $3$-balls and $N_{k}$ is homeomorphic to $S^{2} \times [0,1]$ for each $k \in \{1, \cdots, n-1 \}$. 

We will construct smooth maps $\Phi_{k} : N_{k} \to T_{k}$ as natural extensions of $\psi_{k}$'s by identifying $N_{k-1} \cap N_{k}$ with $F_{k}$ for each type of $T_k$.

\begin{figure}[H]
                \centering
                {\unitlength=1cm
                \begin{picture}(10,14)(0,0.5)
                \put(0,0.5){\includegraphics[width=0.8\textwidth,clip]{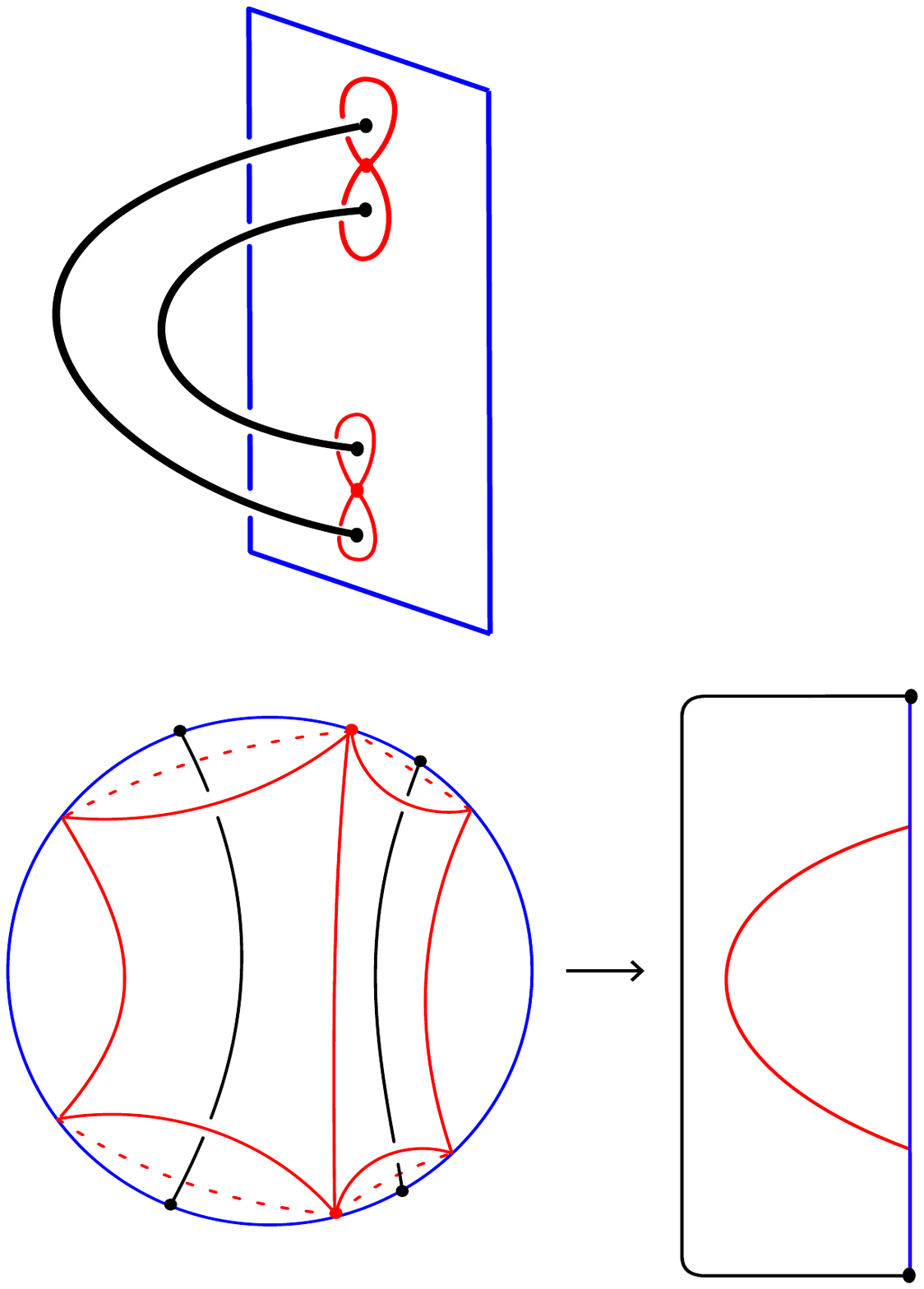}}
                \put(3,7.3){\rotatebox{90}{\scalebox{3}{$\simeq$}}}

                \put(3,14){$F_{1}$}
                \put(4.5,13.2){$1$}
                \put(4.5,12.3){$2$}
                \put(4.4,9.7){$3$}
                \put(4.4,8.9){$4$}

                \put(1,6.8){$N_{0}$}
                \put(2,7){$1$}
                \put(4.6,6.8){$2$}
                \put(1.9,1.2){$3$}
                \put(4.5,1.2){$4$}
                \put(6.5,4.5){$\Phi_{0}$}
                \put(7.7,6.7){$T_{0}$}
                \put(10.2,4){$\gamma_{1}$}
                \put(10.2,7.1){$1,2$}
                \put(10.2,0.9){$3,4$}
                \end{picture}}
                \caption{A smooth map $\Phi_{0}$ from $N_0$ to the region $T_{0}$ of Type $1$.}
                \label{70-4}
            \end{figure}

For $N_0$ and the rectangular region $T_0$ of Type 1, we construct a smooth map $\Phi_{0} : N_{0} \to T_{0}$ as a natural extension of $\psi_{1}$ by identifying $F_1$ with $N_{0} \cap N_{1}$ as depicted in 
Fig.~\ref{70-4}. 

In the figure, the two arcs embedded in $N_{0}$ indicate parts of the components of the two-bridge link $L$, and each point on the arcs is a local singularity of definite fold type. 
The image of these arcs equals the intersection $\partial E \cap T_0$ of the boundary $\partial E$ of the region $E$ and $T_0$. 
The two curves with singularities on $F_1 = \partial N_0$ depict fibers each containing a singular point of indefinite fold type. 

Also, for $N_{n}$ and the rectangular region $T_n$ of Type 4, we can construct a smooth map $\Phi_{n} : N_{n} \to T_{n}$ in the same way.

\begin{figure}[H]
                \setlength\unitlength{1truecm}
                \begin{picture}(12.5,10.5)(0,0.2)
                    \put(0,0.5){\includegraphics[width=1.0\textwidth,clip]{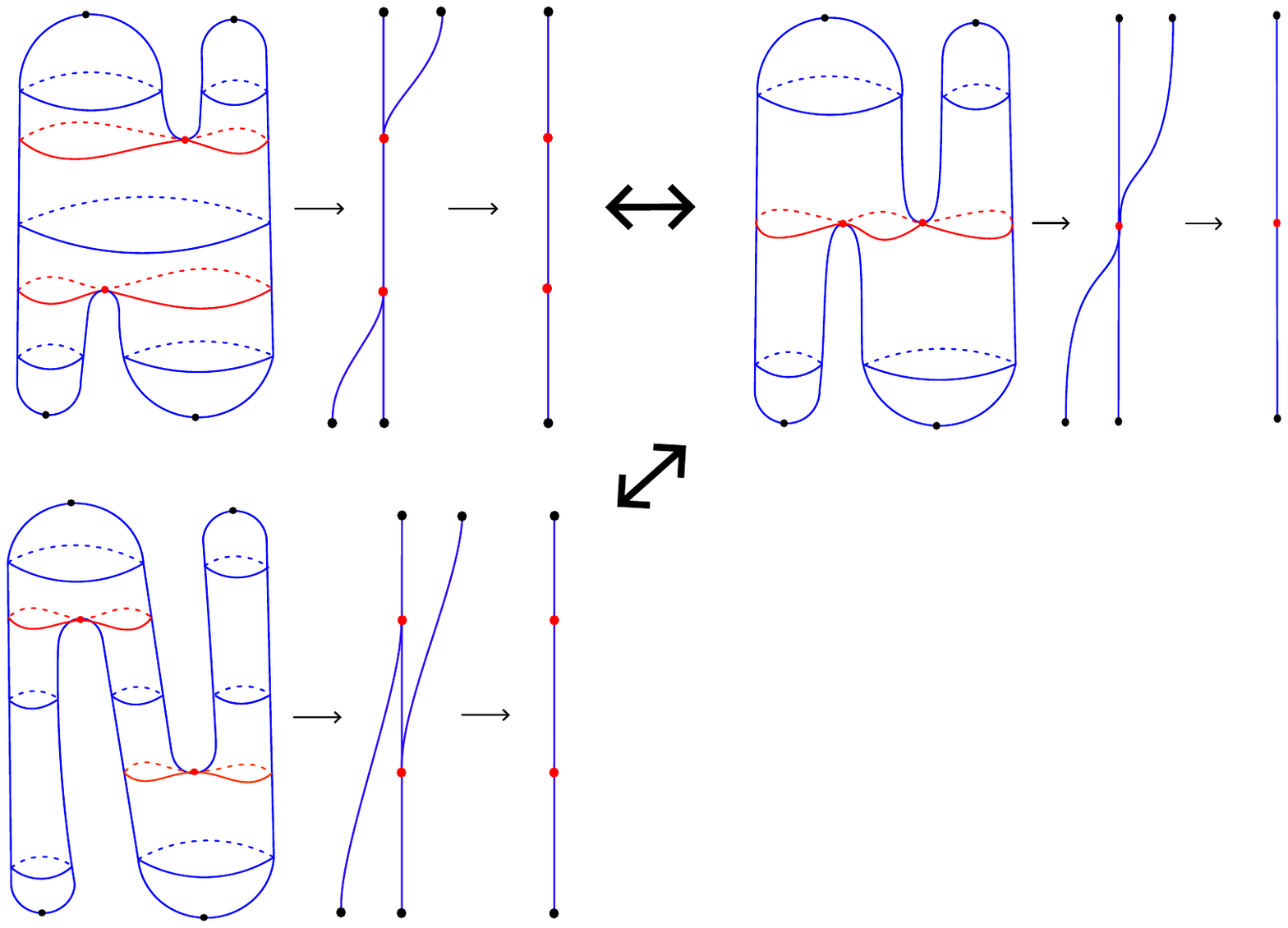}}

                    \put(1,9.2){$1$}
                        \put(2.4,9.2){$2$}
                        \put(0.6,5.1){$3$}
                        \put(2,5.1){$4$}
    
                        \put(3.7,9.2){$1$}
                        \put(4.3,9.2){$2$}
                        \put(3.3,5.1){$3$}
                        \put(3.8,5.1){$4$}
    
                        \put(5.2,9.2){$1,2$}
                        \put(5.2,5.05){$3,4$}
                    \put(7.8,9.2){$1$}
                        \put(9.2,9.2){$2$}
                        \put(7.4,5){$3$}
                        \put(8.8,5){$4$}
    
                        \put(10.4,9.2){$1$}
                        \put(11,9.2){$2$}
                        \put(10,5.1){$3$}
                        \put(10.5,5.1){$4$}
    
                        \put(11.8,9.2){$1,2$}
                        \put(11.8,5.05){$3,4$}
                    \put(1,4.7){$1$}
                        \put(2.4,4.6){$2$}
                        \put(0.6,0.5){$3$}
                        \put(2,0.5){$4$}
    
                        \put(3.7,4.6){$1$}
                        \put(4.3,4.6){$2$}
                        \put(3.3,0.5){$3$}
                        \put(3.8,0.5){$4$}
    
                        \put(5.2,4.6){$1,2$}
                        \put(5.2,0.5){$3,4$}
                \end{picture}
                \caption{A deformation between the Morse functions $\psi_{k}$ and $\psi_{k+1}$.}
                \label{93}
            \end{figure}

\begin{figure}[H]
                \setlength\unitlength{1truecm}
                \begin{picture}(12.5,6)(0,0)
                    \put(-2,-2.5){\includegraphics[width=1.3\textwidth,clip]{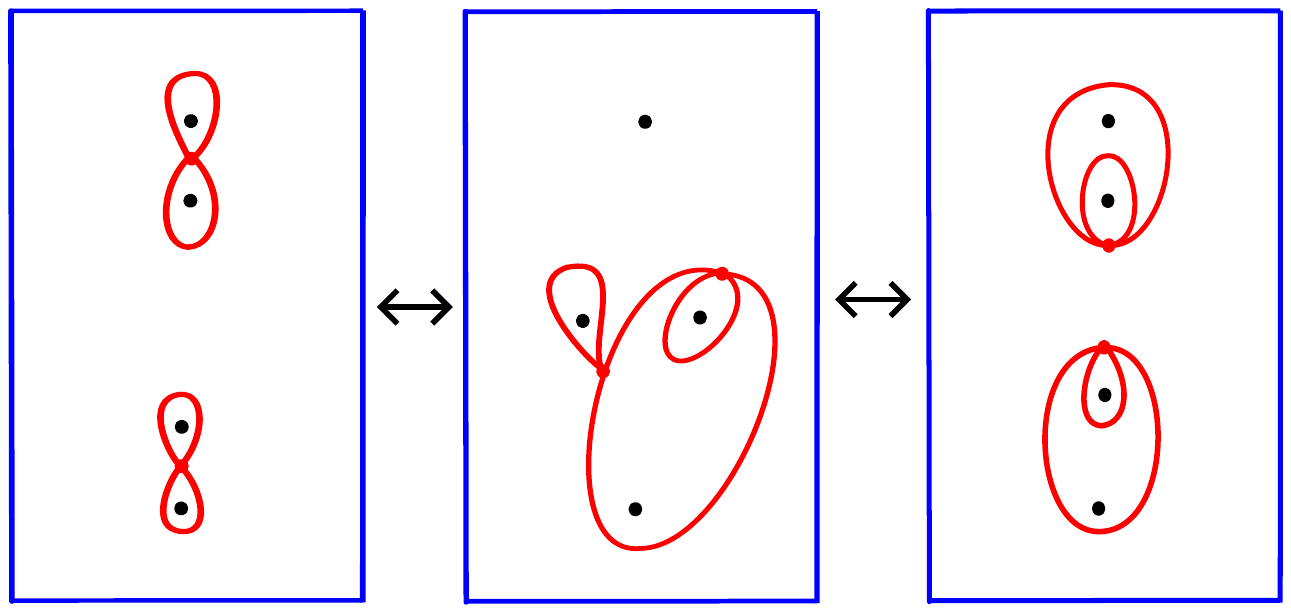}}
                    \put(0.5,5.6){$F_{k} = F_{k+1} $}
                    \put(2.4,5){$1$}
                    \put(2.4,4){$2$}
                    \put(2.2,2){$3$}
                    \put(2.2,1.2){$4$}
                    \put(4.7,5.5){$F_{k}^{\prime}$}
                    \put(6.5,5){$1$}
                    \put(7,2.5){$2$}
                    \put(5.5,2.5){$3$}
                    \put(6.3,1.2){$4$}
                    \put(9,5.5){$F_{k+1}^{\prime \prime}$}
                    \put(10.7,4.9){$1$}
                    \put(10.8,4.5){$3$}
                    \put(10.8,2){$2$}
                    \put(10.6,1.4){$4$}
                \end{picture}
                \caption{Moves of the singular points. 
                }
                \label{931}
            \end{figure}

For $N_k$ and the rectangular region $T_k$ of Type 2, we construct a smooth map $\Phi_{k} : N_{k} \to T_{k}$ as natural extensions of $\psi_{k}$ and $\psi_{k+1}$ by identifying $F_k$ and $F_{k+1}$ with $N_{k-1} \cap N_{k}$ and $N_{k} \cap N_{k+1}$, respectively. 
The natural extension is obtained by a deformation between the Morse functions from $\psi_{k}$ on $F_k \simeq N_{k-1} \cap N_{k}$ to $\psi_{k+1}$ on $F_{k+1} \simeq N_{k} \cap N_{k+1}$. 
The deformation is described in Fig.~\ref{93}. 
It starts at the top left corresponding to $\psi_{k}$, goes right, goes down left, and gets back conversely to the top left corresponding to $\psi_{k+1}$. 
We denote the intermediate sphere (top right) by $F_{k}^{\prime}$ and the intermediate sphere (bottom left) by $F_{k+1}^{\prime \prime}$.
Note that, during the deformation, the singular points and the singular fibers move on the spheres as shown in Fig.~\ref{931}. 
In the figure, moves of the singular points $1,2,3,4$ and the singular fibers on $F_k \simeq N_{k-1} \cap N_{k}$, $F^{\prime}_{k}$, $F_{k+1}^{\prime \prime}$, and $F_{k+1} \simeq N_{k} \cap N_{k+1}$ are illustrated. 

\begin{figure}[H]
                \setlength\unitlength{1truecm}
                \begin{picture}(12.5,6.5)(0,-0.5)
                    \put(0,-1.5){\includegraphics[width=1.0\textwidth,clip]{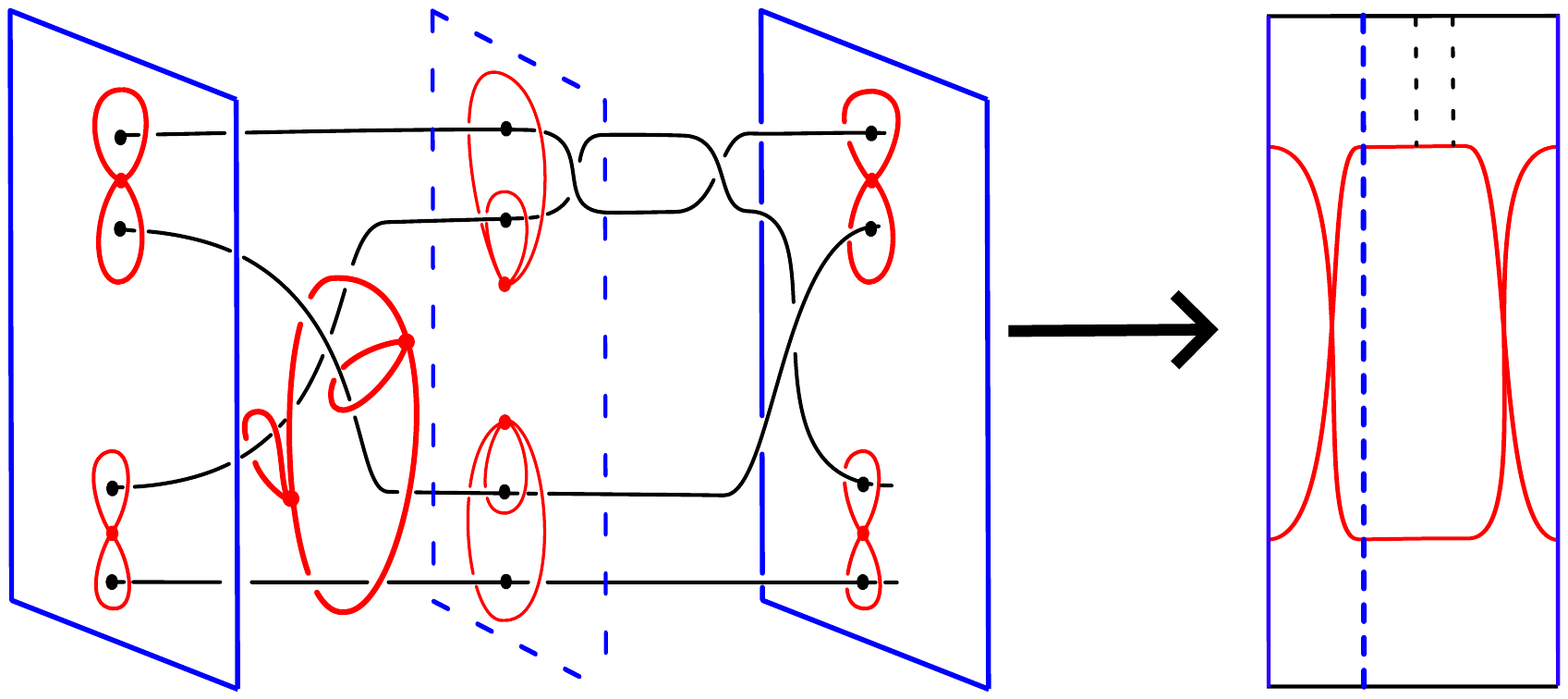}}

                    \put(1.4,4.9){$F_{k}$}
                    \put(0.7,4.2){$1$}
                    \put(0.7,3.1){$2$}
                    \put(0.7,1.8){$3$}
                    \put(0.7,1){$4$}
                        \put(4.2,5){$F_{k+1}^{\prime \prime}$}
                        \put(7,5){$F_{k+1}$}
                    \put(8.5,3.1){$\Phi_{k}$}
                    \put(9.8,5.5){$\gamma_{k}$}
                    \put(10.5,5.5){$\gamma_{k+1}^{\prime \prime}$}
                    \put(12,5.5){$\gamma_{k+1}$}
                    \end{picture}
                \caption{A smooth map $\Phi_{k} : N_{k} \to T_{k}$ for $T_k$ of Type 2.}
                \label{89}
            \end{figure}

Using the deformation given in Fig.~\ref{93} and moves of the points and curves shown in Fig.~\ref{931}, we can illustrate the smooth map $\Phi_{k} : N_{k} \to T_{k}$ in Fig.~\ref{89}. 
Remark that the figure exhibits the map in the case of $|b_{i}|=2$. 

In the figure, the four arcs connecting $F_k$ and $F_{k+1}$ embedded in $N_{k}$ are parts of the two-bridge link $L$, and each point on the arcs is a local singularity of definite fold type. 
The image of these arcs equals the intersection $\partial E \cap T_k$. 
Note that, for the smooth map $\Phi_{k} : N_{k} \to T_{k}$, there exist exactly two fibers of type $\mathrm{I\hspace{-1.2pt}I^{2}}$ in $N_{k}$, each of which contains two indefinite fold points. 

\begin{figure}[H]
                    \centering
                    {\unitlength=1cm
                    \begin{picture}(12,7.5)(0,-0.5)
                    \put(1,0){\includegraphics[width=0.7\textwidth,clip]{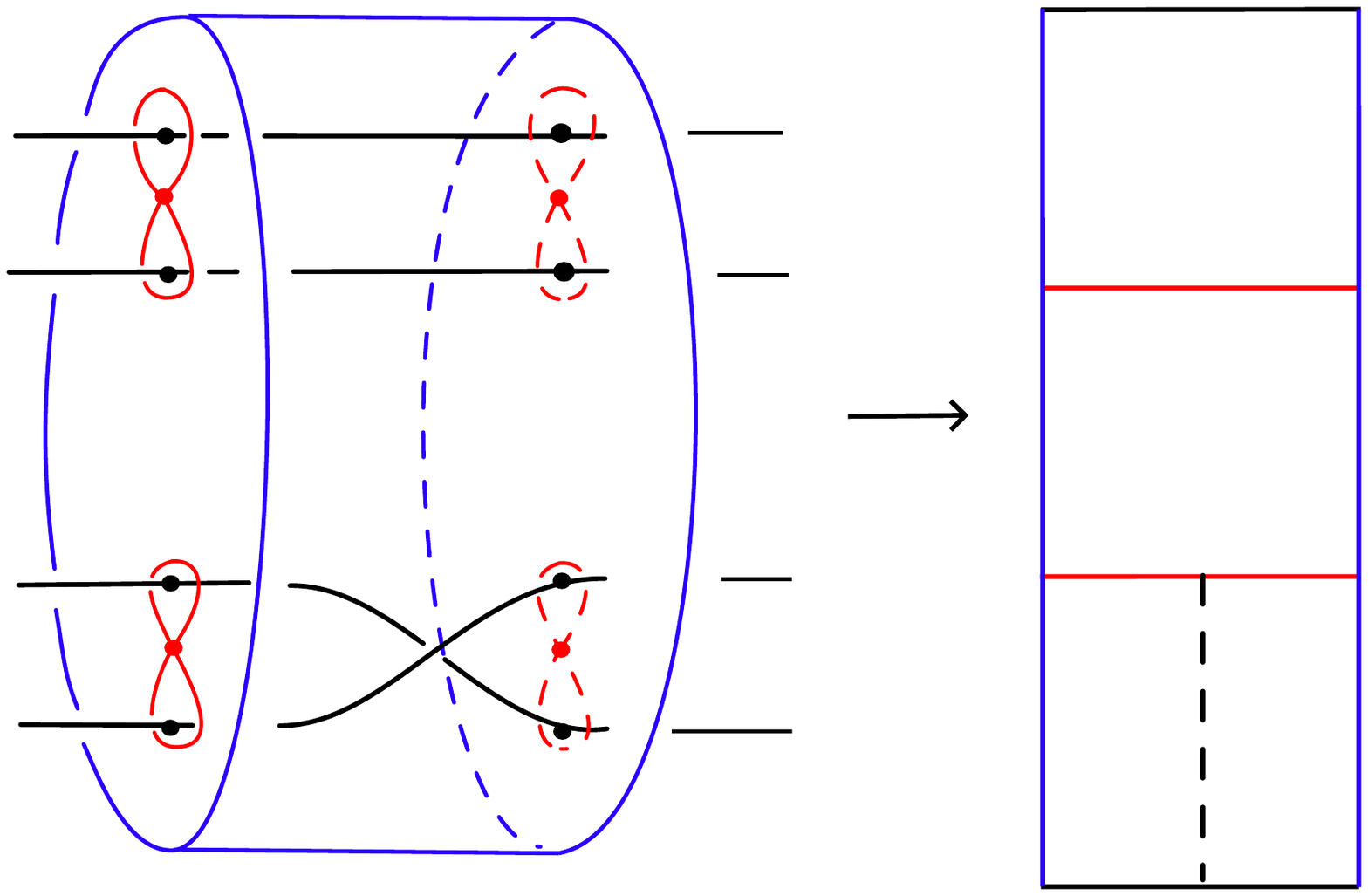}}
                    \put(1,6.3){\underline{\textbf{Type $3$}}}
                    \put(6.65,3.5){$\Phi_{k}$}
                    \put(1.3,5.5){$N_{k}$}
                    \put(7,5.5){$T_{k}$}

                    \put(1,5){$1$}
                    \put(1,4){$2$}
                    \put(1,2.2){$3$}
                    \put(1,1.3){$4$}

                    \put(2,0.1){$F_{k}$}
                    \put(4.8,0.1){$F_{k+1}$}
                    \put(7.5,0){$\gamma_{k}$}
                    \put(9.5,0){$\gamma_{k+1}$}
                    \end{picture}}
                    \caption{ A smooth map to the region of Type $3$. }
                    \label{type3}
                \end{figure}

For $N_k$ and the rectangular region $T_k$ of Type 3, in a similar way as before, we construct a smooth map $\Phi_{k} : N_{k} \to T_{k}$ as natural extensions of $\psi_{k}$ and $\psi_{k+1}$ by identifying $F_k$ and $F_{k+1}$ with $N_{k-1} \cap N_{k}$ and $N_{k} \cap N_{k+1}$, respectively. 
See Fig.~\ref{type3}. 
In this case, we note that there are no singular fibers containing two indefinite fold points.

Now, we connect $N_{0}, N_{1}, \cdots, N_{n-1}$ and $N_{n}$ by identifying the boundaries of them to obtain $S^{3}$. 
Also, we connect $T_{0}, T_{1}, \cdots, T_{n-1}$ and $T_{n}$ by identifying $\gamma_{1}, \cdots, \gamma_{n}$ to make $E$. 
Then, a smooth map $f_{2} : S^{3} \to \mathbb{R}^{2}$ satisfying $f_{2} |_{N_{k}} = \Phi_{k}$ for $k \in \{ 0,1,2,\cdots, n \}$ is constructed. 
By construction, this map $f_2$ is a stable map from $S^3$ into $\mathbb{R}^{2}$ with desired properties. 
In particular, $f_2$ has no singular fibers of type $\mathrm{I\hspace{-1.2pt}I^{3}}$ and the total number of singular fibers of type $\mathrm{I\hspace{-1.2pt}I^{2}}$ is $2m$. 
This completes a proof of Theorem~\ref{main1}.

    \begin{example}
        A stable map $f : S^{3} \to \mathbb{R}^2$ with the figure-eight knot as $S_{0}(f)$ obtained by the construction given in this section is presented in Fig.~\ref{94}. 
        This $f$ has just two fibers of type $\mathrm{I\hspace{-1.2pt}I^{2}}$, each of which contains two indefinite fold points. 
        Some of the other fibers containing an indefinite fold point are also shown in the figure. 

        \begin{figure}[H]
                \setlength\unitlength{1truecm}
                \begin{picture}(12.5,14)(0,0)
                    \put(2,0){\includegraphics[width=0.7\textwidth,clip]{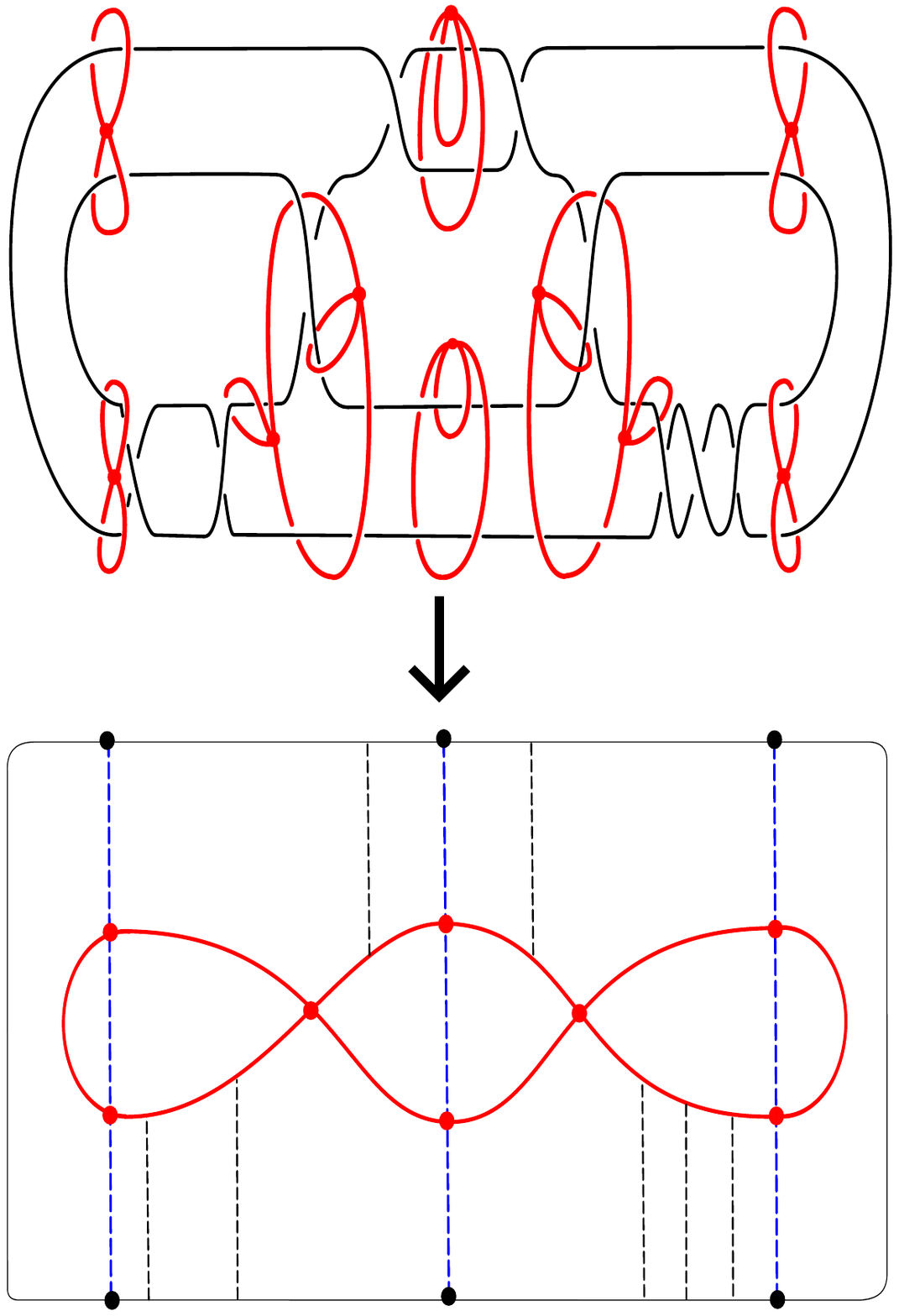}}
                \end{picture}
                \caption{A stable map $f$ from $S^{3}$ with the figure-eight knot as $S_{0}(f)$.}
                \label{94}
        \end{figure}

    \end{example}

\section{A construction of $f_{3}$}

In this section, we give a construction of 
a stable map $f_{3}$ from $S^{3}$ into $\mathbb{R}^{2}$ satisfying that $S_{0}(f_{3})= L$, $\mathrm{I\hspace{-1.2pt}I^{2}}(f_{3}) = \emptyset$ and $|\mathrm{I\hspace{-1.2pt}I^{3}}(f_{3})| = \frac{1}{2} \sum_{i=1}^{m} | b_{i} |$
for a two-bridge link $L$ in $S^3$ when $L$ has a Conway form $C(a_{1},b_{1},\cdots,a_{m},b_{m},a_{m+1} )$ with all the $b_{i}$'s are even. 

Let $D^{\prime}=  C(a_{1},b_{1},\cdots,a_{m},b_{m},a_{m+1} )$ be a Conway form of $L$ shown in Fig.~\ref{82} assuming that $b_{i}$ is an even number for all $i$. 
As in the previous section, we set a rectangular region $E$ on the plane $\mathbb{R}^{2}$ which contains $D^{\prime}$ inside. 

\begin{figure}[H]
                \setlength\unitlength{1truecm}
                \begin{picture}(12.5,4.5)(0,0.3)
                    \put(1,-1){\includegraphics[height=7cm,clip]{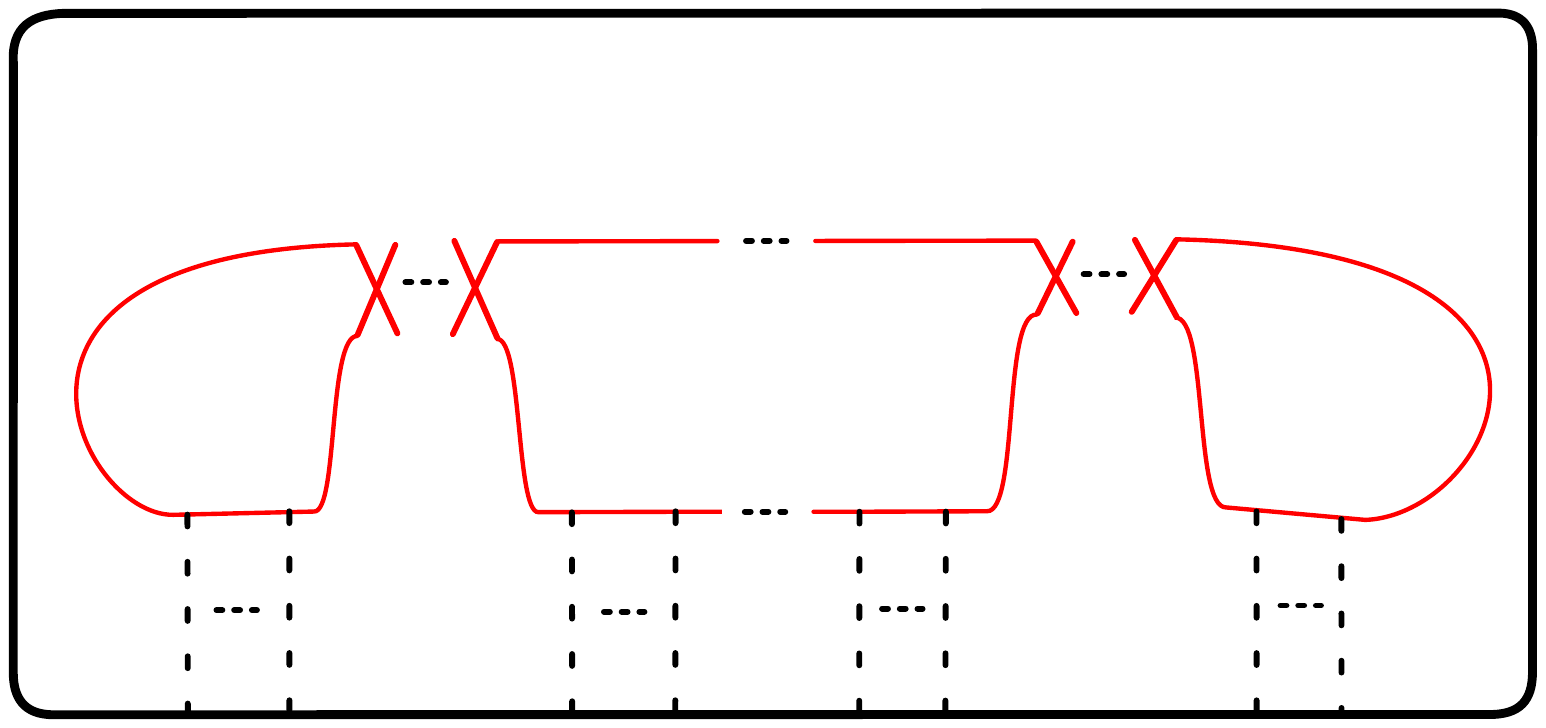}}
                    \put(3.65,3.7){$\overbrace{ }^{|b_{1}|}$}
                    \put(7.5,3.7){$\overbrace{ }^{|b_{m}|}$}

                    \put(2,5){$E$}
                    \put(2,3.4){$\hat{C^{\prime}}$}

                    \put(2.6,0.75){$\underbrace{ }_{|a_{1}|}$}
                    \put(4.77,0.75){$\underbrace{ }_{|a_{2}|}$}
                    \put(6.35,0.75){$\underbrace{ }_{|a_{m}|}$}
                    \put(8.47,0.75){$\underbrace{ }_{|a_{m+1}|}$}
                \end{picture}
                \caption{The immersed curve $\hat{C^{\prime}}$.}
                \label{1}
            \end{figure}

Let $\hat{C^{\prime}}$ be an immersed curve obtained from $D^{\prime}$ as in the previous section: 
Perform horizontal smoothing at all the double points adjacent to the outer region, removing the outermost circle appearing after the smoothing, and forgetting the crossing information of the remaining diagram. 
See Fig.~\ref{1} similar as Fig.\ref{90}. 
Note that the immersed curve $\hat{C^{\prime}}$ has $\sum_{i=1}^{m} | b_{i} |$ normal crossings. 

\begin{figure}[H]
                \setlength\unitlength{1truecm}
                \begin{picture}(12.5,2.8)(0,1)
                    \put(1.5,-1){\includegraphics[width=.75\textwidth,clip]{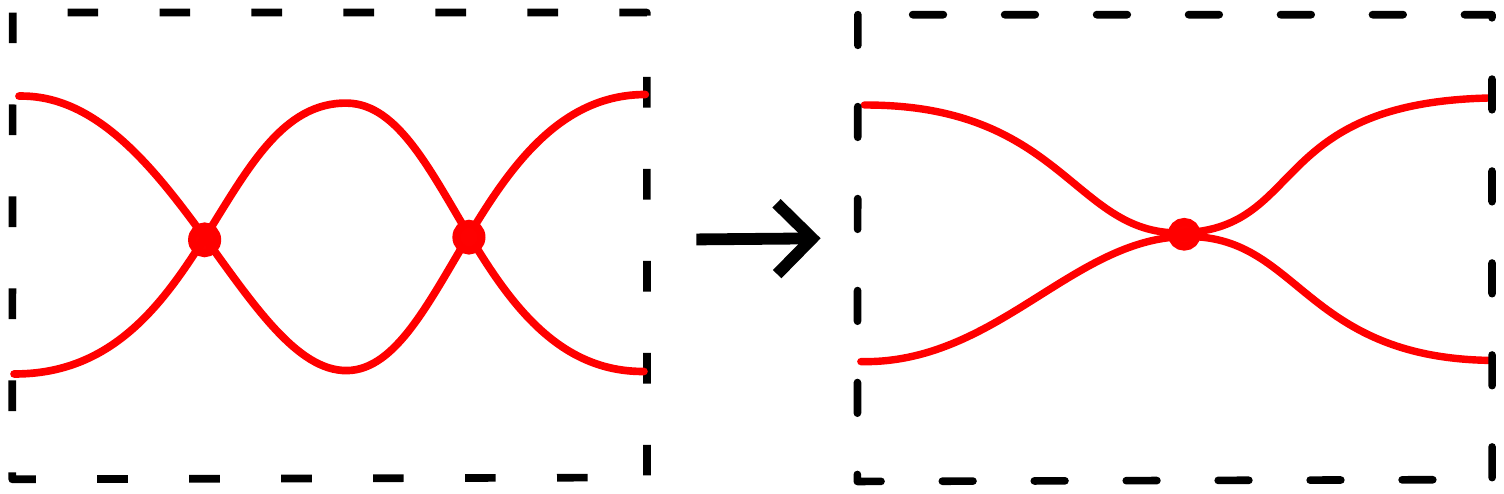}}
                \end{picture}
                \caption{Modification on the regions.}
                \label{replace}
            \end{figure}

Then, we perform modifications on the immersed curve $\hat{C^{\prime}}$ at bigon regions as illustrated in Fig.~\ref{replace}
so that the obtained curve $C^{\prime}$ has no normal crossings and has $ \frac{1}{2} \sum_{i=1}^{m} | b_{i} |$ self tangent points. 
Remark that this can be done since all the $b_{i}$'s are assumed to be even. 

\begin{figure}[H]
                \setlength\unitlength{1truecm}
                \begin{picture}(12.5,5)(0,0)
                    \put(2.5,0){\includegraphics[height=5cm,clip]{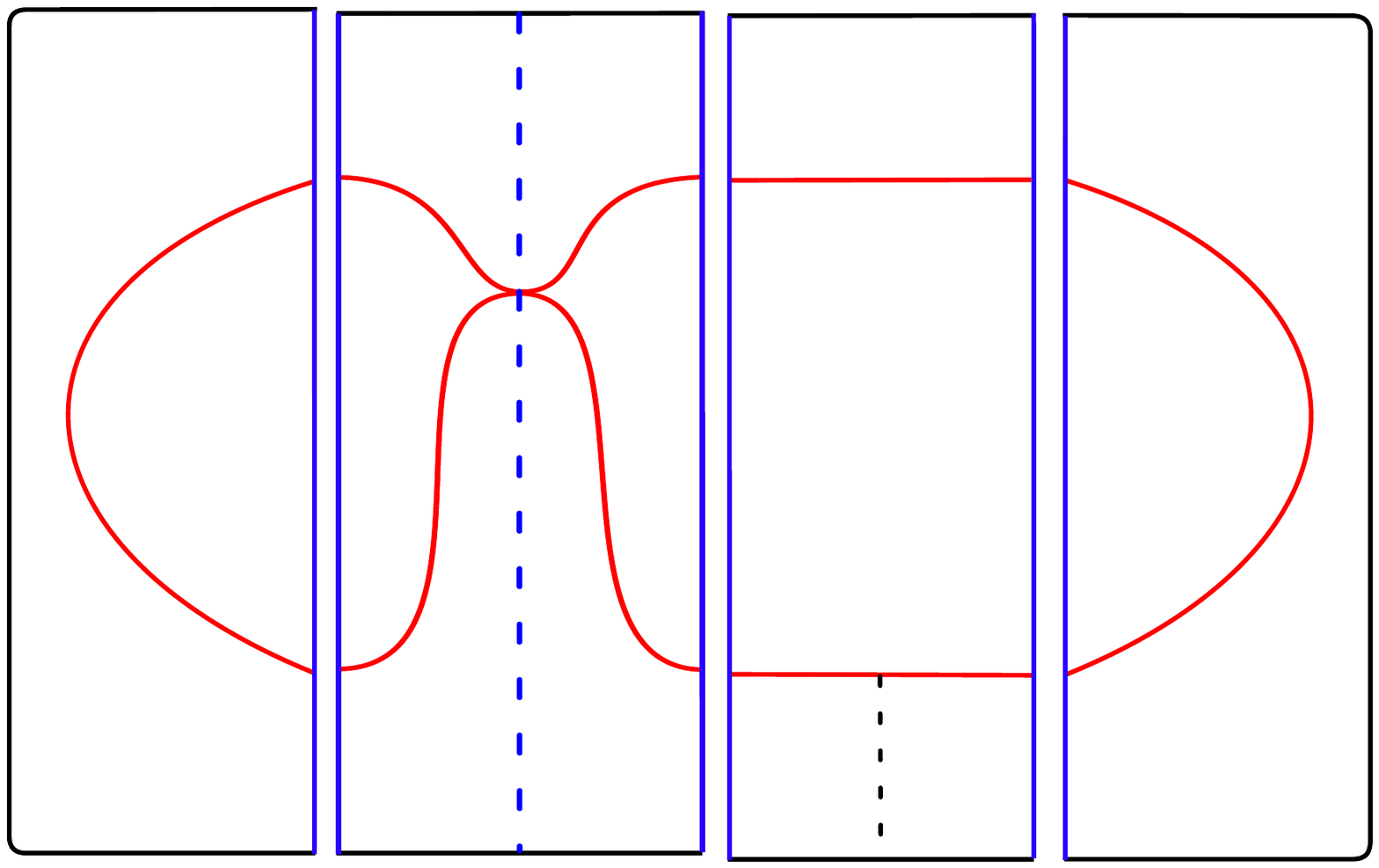}}
                    \put(2.9,5){Type $1$}
                    \put(4.6,5){Type $2$}
                    \put(6.5,5){Type $3$}
                    \put(8.2,5){Type $4$}
                \end{picture}
                \caption{The four types of rectangular regions.}
                \label{2}
            \end{figure}

We take parallel segments $\gamma_{1},\gamma_{2},\cdots,\gamma_{n}$ on $E$ such that they separate $E$ into strips, say, $T_{0},T_{1},T_{2},\cdots,T_{n}$, so that $T_{k} \cap C^{\prime}$ looks like one of the shapes depicted in Fig.~\ref{2} for each $k \in \{0,1,\cdots,n \}$. 
We refer to the rectangular regions separated by $\gamma_{k}$ as Type $1$, Type $2$, Type $3$ and Type $4$ as shown in Fig.~\ref{2}. 
In the regions of Type 3, dotted lines are added to mark the crossings in the twists on $D^{\prime}$ associated to the $a_i$'s. 
The region of Type 2 only contains a self tangent point of $C^{\prime}$. 
Thus, the number of the regions of Type 2 is equal to $\frac{1}{2} \sum_{i=1}^{m} | b_{i} |$. 

\begin{figure}[H]
                \centering
                 {\unitlength=1cm
                \begin{picture}(10,7)(0,0.5)
                \put(0,0.5){\includegraphics[width=0.8\textwidth,clip]{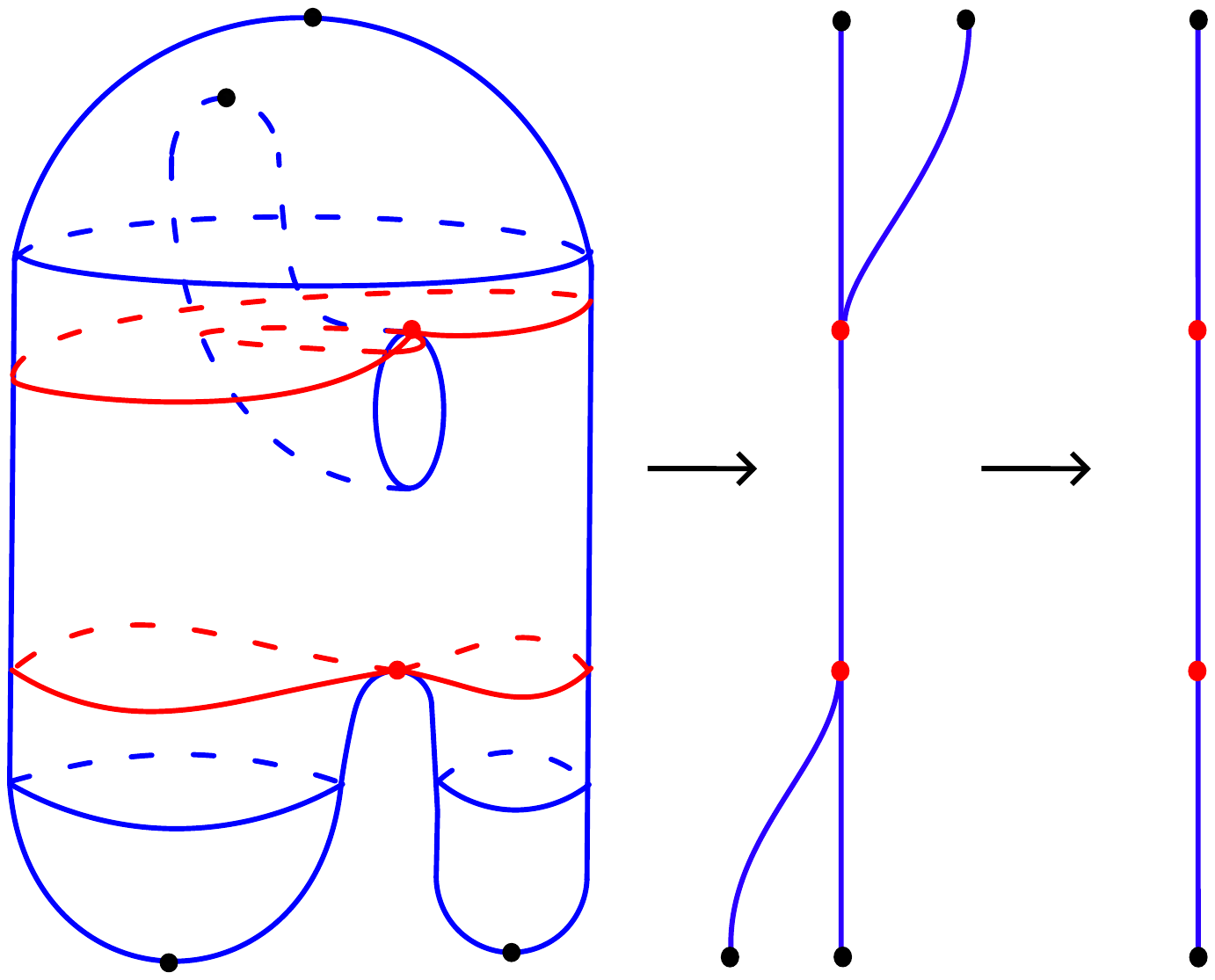}}
                \put(0,4){$\psi_{k} :$}
                \put(0.5,6.5){$F_{k}$}
                \put(5.5,6.5){$\tau_{k}$}
                \put(8,6.5){$\gamma_{k}$}

                \put(3,7.3){$1$}
                \put(2.5,6.5){$2$}
                \put(1.7,0.5){$3$}
                \put(4,0.5){$4$}

                \put(6.1,7.3){$1$}
                \put(7,7.3){$2$}
                \put(5.5,0.5){$3$}
                \put(6.1,0.5){$4$}

                \put(8.3,7.3){$1,2$}
                \put(8.3,0.5){$3,4$}
                \end{picture}}
                \caption{ A Morse function $\psi_{k} : F_{k} \to \gamma_{k} \subset E$.}
                \label{4}
            \end{figure}

For each $k \in \{ 1,2,\cdots,n \}$, we construct a smooth map $\psi_{k} : F_{k} \to \gamma_{k}$ whose Reeb graph is the train track $\tau_{k}$ as shown in Fig.~\ref{4}, where $F_{k}$ is homeomorphic to the 2-sphere $S^{2}$. 
In the figure, on $F_{k}$, the circles show regular fibers and the singular curves show singular fibers.
Note that $\psi_{k}$ gives a Morse function on $F_{k} \simeq S^{2}$. 

\begin{figure}[H]
                \centering
                {\unitlength=1cm
                \begin{picture}(10,10.5)(0,0)
                \put(-2,0.5){\includegraphics[width=1.1\textwidth,clip]{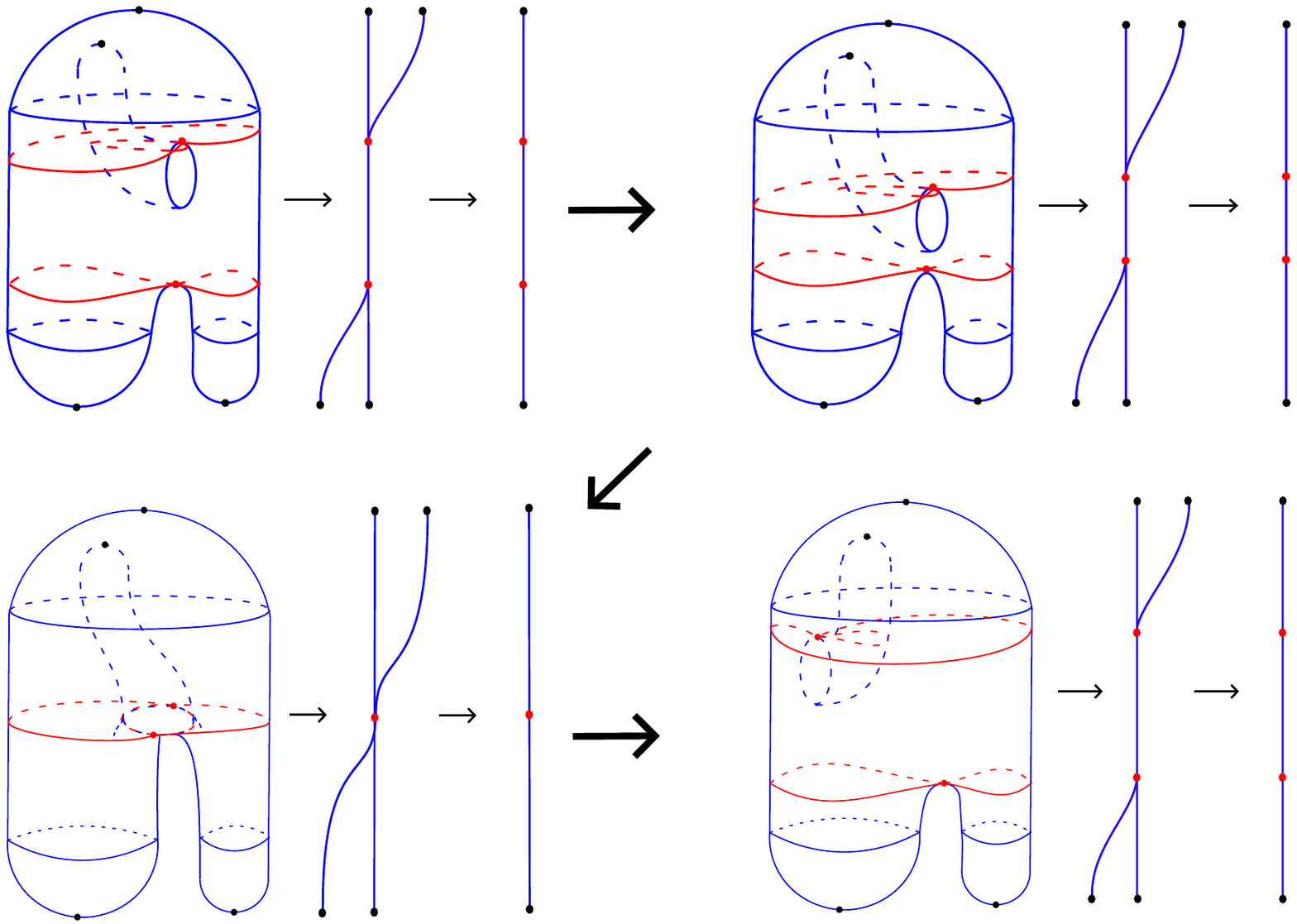}}

                \put(-0.4,10.2){$1$}
                \put(-0.5,9.6){$2$}
                \put(-1,5.6){$3$}
                \put(0.5,5.6){$4$}

                \put(1.9,10.3){$1$}
                \put(2.5,10.3){$2$}
                \put(1.4,5.6){$3$}
                \put(1.9,5.6){$4$}

                \put(3.4,10.3){$1,2$}
                \put(3.4,5.6){$3,4$}

                \put(7.5,10.3){$1$}
                \put(7.3,9.6){$2$}
                \put(6.7,5.6){$3$}
                \put(8.4,5.6){$4$}

                \put(9.8,10.3){$1$}
                \put(10.4,10.3){$2$}
                \put(9.3,5.6){$3$}
                \put(9.9,5.6){$4$}

                \put(11.4,10.3){$1,2$}
                \put(11.4,5.6){$3,4$}

                \put(-0.3,5){$1$}
                \put(-0.5,4.4){$2$}
                \put(-0.9,0.3){$3$}
                \put(0.6,0.3){$4$}

                \put(2,5){$1$}
                \put(2.6,5){$2$}
                \put(1.4,0.3){$3$}
                \put(2,0.3){$4$}

                \put(3.5,5){$1,2$}
                \put(3.5,0.3){$3,4$}

                \put(7.7,5){$1$}
                \put(7.5,4.5){$2$}
                \put(6.8,0.3){$3$}
                \put(8.6,0.3){$4$}

                \put(10,5.1){$1$}
                \put(10.5,5.1){$2$}
                \put(9.5,0.3){$3$}
                \put(10,0.3){$4$}

                \put(11.3,5.1){$1,2$}
                \put(11.3,0.3){$3,4$}
                \end{picture}}
                \caption{A deformation 
                from $\psi_{k}$ to $\psi_{k+1}$. }
                \label{type21}
            \end{figure}

Actually, the Morse function $\psi_{k}$ is the same as that given in the previous section shown in Fig.~\ref{70-3}. 
However, by seeing this figure, there is an advantage to describing the deformation of the Morse functions given later (Fig.~\ref{type21}). 

For a standard projection $\mathbb{R}^{3} \to \mathbb{R}^{2}$ to obtain $D^{\prime}$ from $L$, the preimage of a line on $\mathbb{R}^{2}$ containing $\gamma_{k}$ gives an embedded sphere in $S^{3}$. 
We identify it with $F_{k}$ for each $k  \in \{1,2,\cdots,n \}$, and decompose $S^{3}$ into $N_{0},N_{1},N_{2},\cdots,N_{n}$ by the spheres $F_{1},\cdots,F_{n}$. 
Here $N_{0}$ and $N_{n}$ are homeomorphic to $3$-balls and $N_{k}$ is homeomorphic to $S^{2} \times [0,1]$ for each $k \in \{1, \cdots, n-1 \}$. 

Based on a similar strategy as in the previous section, we will construct smooth maps $\Phi_{k} : N_{k} \to T_{k}$ as natural extensions of $\psi_{k}$'s by identifying $N_{k-1} \cap N_{k}$ with $F_{k}$ for each type.

\begin{figure}[H]
                \centering
                {\unitlength=1cm
                \begin{picture}(10,5)(0,0.5)
                \put(0,0){\includegraphics[width=0.75\textwidth,clip]{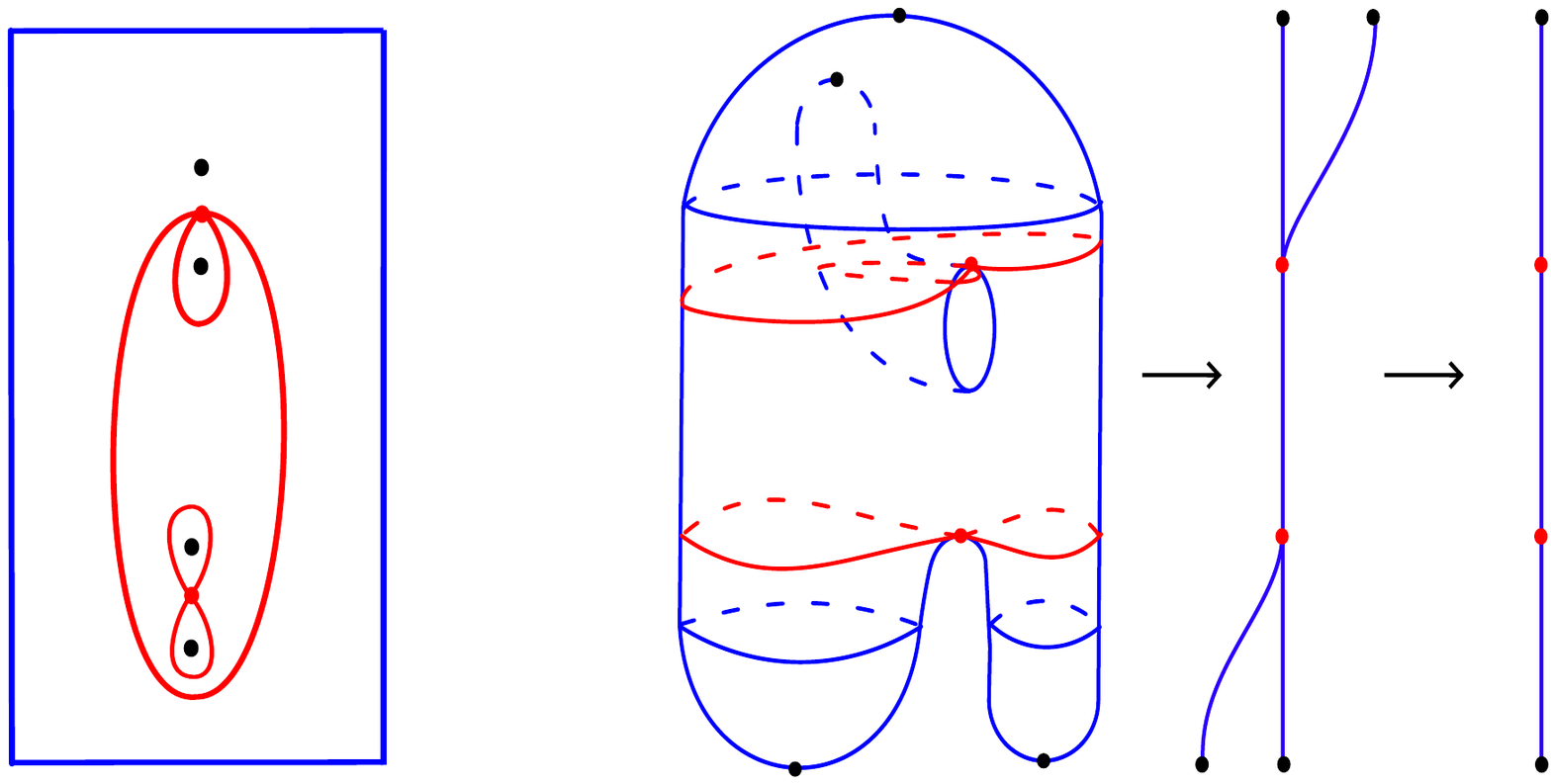}}
                \put(3,3){\scalebox{3}{$\simeq$}}

                \put(1.35,4.7){$1$}
                \put(1.4,3.65){$2$}
                \put(1.35,2.45){$3$}
                \put(1.35,1.85){$4$}

                \put(5.6,5.6){$1$}
                \put(5.2,5.2){$2$}
                \put(4.8,0.6){$3$}
                \put(6.3,0.6){$4$}

                \put(7.75,5.65){$1$}
                \put(8.35,5.65){$2$}
                \put(7.25,0.6){$3$}
                \put(7.8,0.6){$4$}

                \put(9.2,5.7){$1,2$}
                \put(9.2,0.6){$3,4$}
                \end{picture}}
                \caption{The Morse function $\psi_{1}$ on $F_1 \simeq N_{0} \cap N_{1} \subset S^3$. }
                \label{type11}
            \end{figure}

\begin{figure}[H]
                \centering
                {\unitlength=1cm
                \begin{picture}(10,4.5)(0,1)
                \put(0,0){\includegraphics[width=0.8\textwidth,clip]{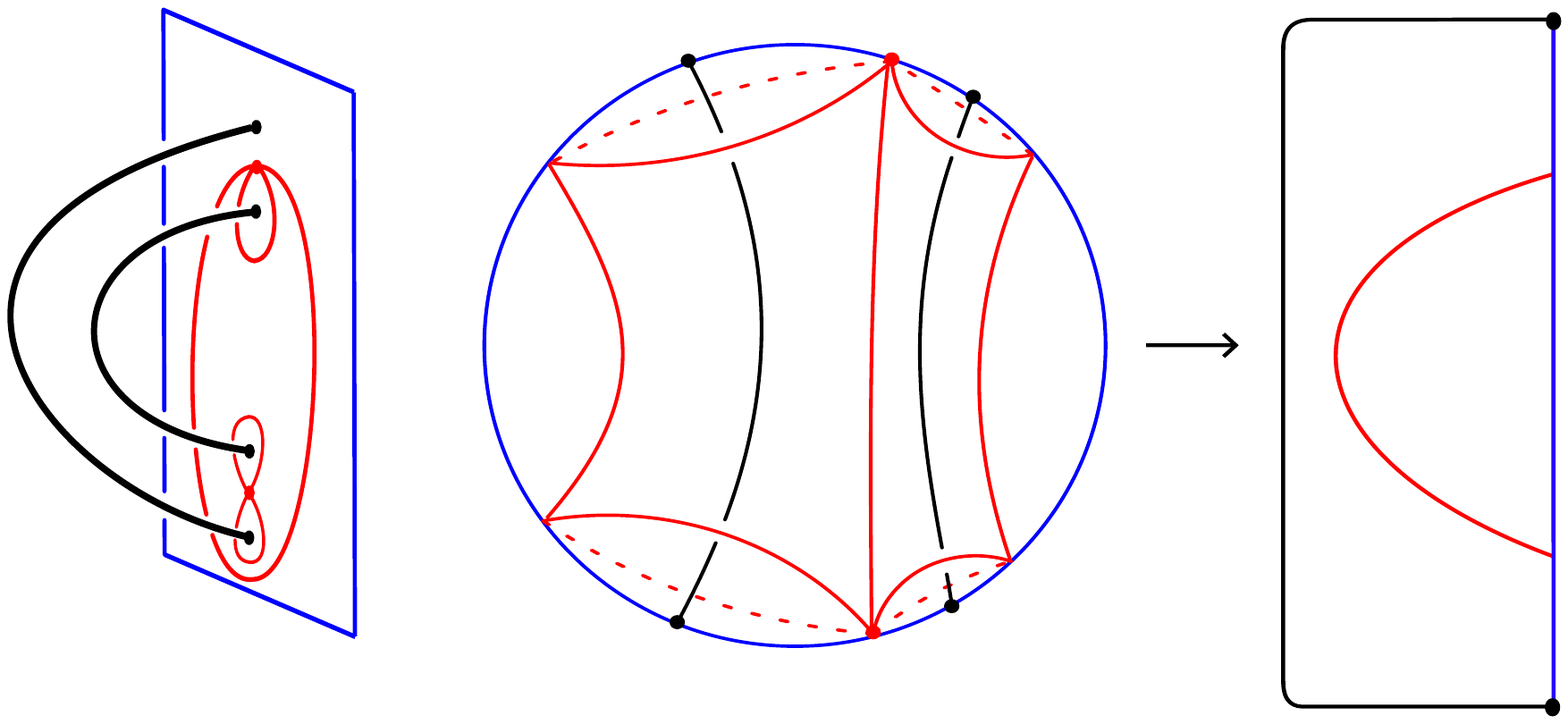}}
                \put(2.45,3.6){\scalebox{1.75}{$\simeq$}}

                \put(2.3,5.3){$F_{1}$}
                \put(1.75,5){$1$}
                \put(1.7,3.9){$2$}
                \put(1.7,3.1){$3$}
                \put(1.7,2.5){$4$}

                \put(3.4,5.5){$N_{0}$}
                \put(4.3,5.6){$1$}
                \put(6.2,5.4){$2$}
                \put(6.2,1.5){$3$}
                \put(4,1.5){$4$}

                \put(8.5,5.3){$T_{0}$}
                \put(8.2,1){$\gamma_{1}$}
                \put(9.8,5.9){$1,2$}
                \put(9.8,0.8){$3,4$}
                \end{picture}}
                \caption{A smooth map from $N_0$ to the region $T_{0}$ of Type $1$. 
                }
                \label{type12}
            \end{figure}

For $N_0$ and the rectangular region $T_0$ of Type 1, we construct a smooth map $\Phi_{0} : N_{0} \to T_{0}$ as a natural extension of $\psi_{1}$ by identifying $F_1$ with $N_{0} \cap N_{1}$ as depicted in Fig.~\ref{type11} and \ref{type12}. 

            
The obtained map $\Phi_{0}$ is the same as that given in the previous section shown in Fig.~\ref{70-4}, but we include Fig.~\ref{type11} and Fig.~\ref{type12} to make easier to see the deformation given later. 
Also, for that reason, $L \cap F_{1}$ are numbered by $1,2,3$ and $4$ in the figures. 

Also, for $N_{n}$ and the rectangular region $T_n$ of Type 4, we can construct a smooth map $\Phi_{n} : N_{n} \to T_{n}$ in the same way.

For $N_k$ and the rectangular region $T_k$ of Type 2, we construct a smooth map $\Phi_{k} : N_{k} \to T_{k}$ as natural extensions of $\psi_{k}$ and $\psi_{k+1}$ by identifying $F_k$ and $F_{k+1}$ with $N_{k-1} \cap N_{k}$ and $N_{k} \cap N_{k+1}$, respectively. 
The natural extension is obtained by a deformation between the Morse functions from $\psi_{k}$ on $F_k \simeq N_{k-1} \cap N_{k}$ to $\psi_{k+1}$ on $F_{k+1} \simeq N_{k} \cap N_{k+1}$. 
In Fig.~\ref{type21},  we describe the deformation between the Morse functions from $\psi_{k}$ to $\psi_{k+1}$.
It starts at the top left corresponding to $\psi_{k}$, goes right, goes down left, and goes right corresponding to $\psi_{k+1}$. 
In Fig.~\ref{type21}, we denote the intermediate spheres (top right and bottom left) by $F^{\prime}_{k}$ and $F_{k+1}^{\prime \prime}$. 

\begin{figure}[H]
                \centering
                {\unitlength=1cm
                \begin{picture}(10,5)(0,0)
                \put(-1.3,-2.2){\includegraphics[width=1\textwidth,clip]{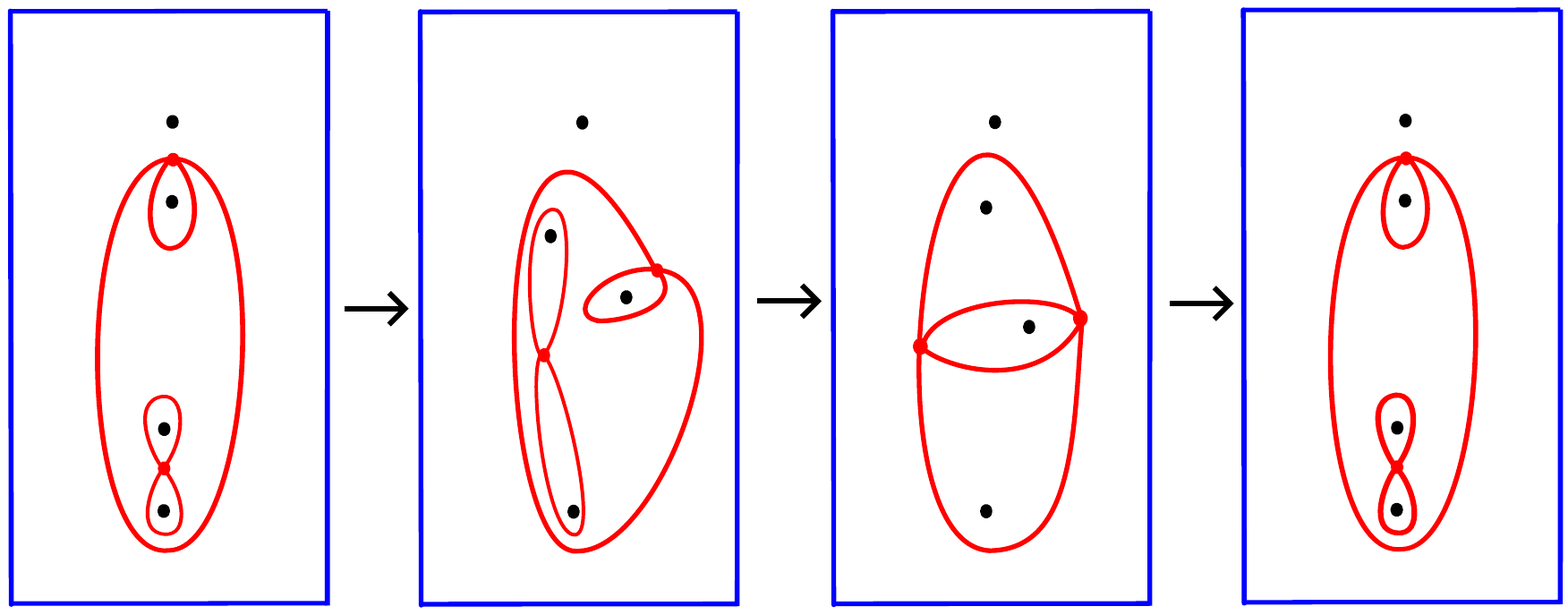}}

                \put(-1,4.4){$F_{k}$}
                \put(0.4,3.8){$1$}
                \put(0.4,2.9){$2$}
                \put(0.3,1.4){$3$}
                \put(0.2,0.8){$4$}

                \put(2.2,4.4){$F^{\prime}_{k}$}
                \put(3.7,3.7){$1$}
                \put(3.8,1.9){$2$}
                \put(2.5,2.6){$3$}
                \put(3.5,1){$4$}

                \put(5.5,4.4){$F_{k+1}^{\prime \prime}$}
                \put(7,3.7){$1$}
                \put(6.7,2){$2$}
                \put(6.8,3){$3$}
                \put(6.8,0.7){$4$}

                \put(9,4.4){$F_{k+1}$}
                \put(10.2,3.7){$1$}
                \put(10.2,2.7){$2$}
                \put(10.2,1.5){$3$}
                \put(10.2,0.8){$4$}
                \end{picture}}
                \caption{Moves of the singular points and the singular fibers on $F_k \simeq N_{k-1} \cap N_{k}$, $F^{\prime}_{k}$, $F_{k+1}^{\prime \prime}$ and $F_{k+1} \simeq N_{k} \cap N_{k+1}$. }
                \label{type22}
            \end{figure}

Note that, during the deformation, the singular points and the singular fibers move on the spheres as shown in Fig.~\ref{type22}.
In the figure, moves of the singular points $1,2,3,4$ and the singular fibers on $F_k \simeq N_{k-1} \cap N_{k}$, $F^{\prime}_{k}$, $F_{k+1}^{\prime \prime}$, and $F_{k+1} \simeq N_{k} \cap N_{k+1}$ are illustrated. 

\begin{figure}[H]
                \centering
                {\unitlength=1cm
                \begin{picture}(10,5.5)(0,0)
                \put(-1.3,-1.5){\includegraphics[width=1\textwidth,clip]{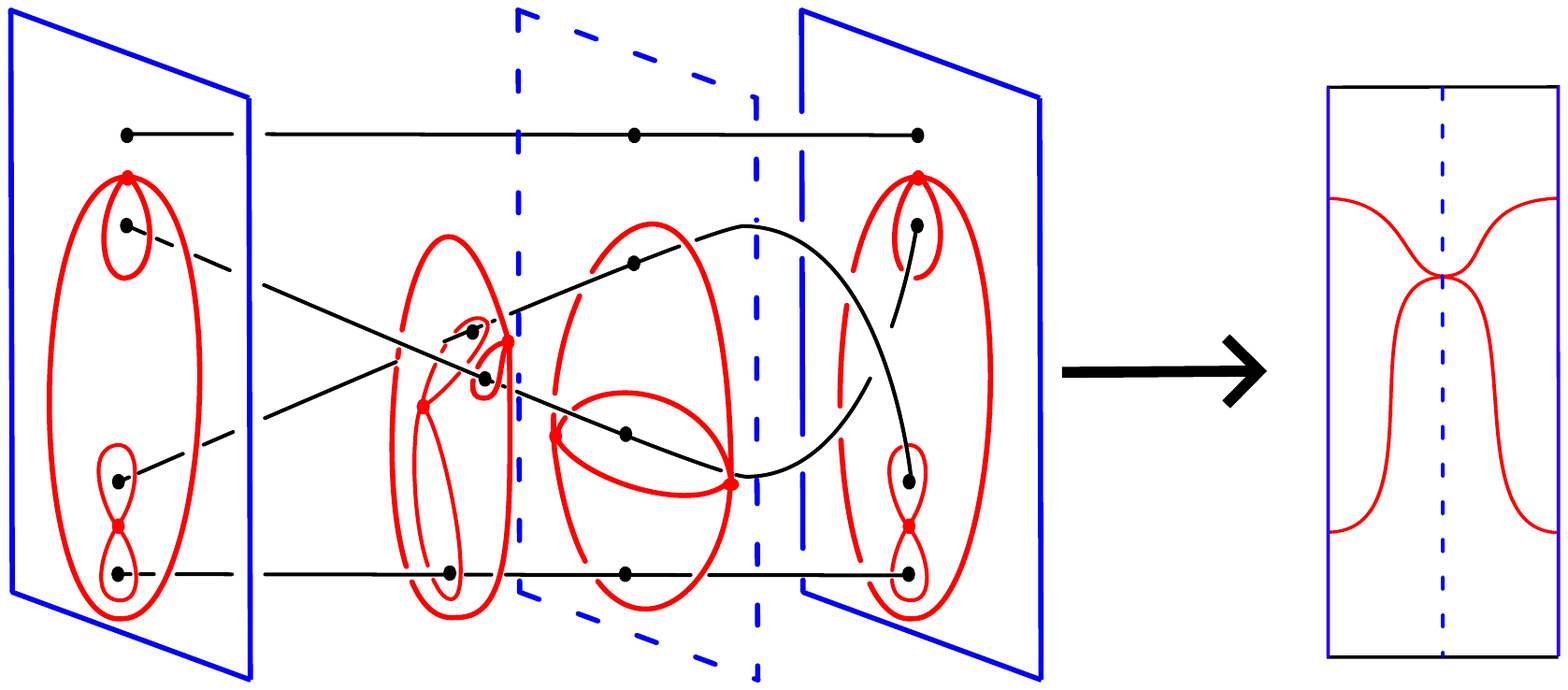}}

                \put(-0.7,4.5){$1$}
                \put(-0.7,3.5){$2$}
                \put(-0.7,2){$3$}
                \put(-0.7,1.3){$4$}

                \put(-1.1,5){$F_{k}$}
                \put(2.9,5){$F_{k+1}^{\prime \prime}$}
                \put(5.1,5){$F_{k+1}$}
                \put(7.5,3.1){$\Phi_{k}$}
                \put(8.9,5.2){$\gamma_{k}$}
                \put(10.5,5.2){$\gamma_{k+1}$}
                \end{picture}}
                \caption{A smooth map $\Phi_{k} : N_{k} \to T_{k}$ for $T_k$ of Type 2. }
                \label{type23}
            \end{figure}

Using the deformation given in Fig.~\ref{type21} and moves of the points and curves shown in Fig.~\ref{type22}, we can illustrate a smooth map $\Phi_{k} : N_{k} \to T_{k}$ as in Fig.~\ref{type23}.

In the figure, the four arcs connecting $F_k$ and $F_{k+1}$ embedded in $N_{k}$ are parts of the two-bridge link $L$, and each point on the arcs is a local singularity of definite fold type. 
The image of these arcs equals the intersection $\partial E \cap T_k$. 

Note that, for the smooth map $\Phi_{k} : N_{k} \to T_{k}$, there exists exactly one singular fiber of type $\mathrm{I\hspace{-1.2pt}I^{3}}$ lying on $F_{k+1}^{\prime \prime}$ in $N_{k}$ and there are no singular fibers of type $\mathrm{I\hspace{-1.2pt}I^{2}}$. 

\begin{figure}[H]
                    \centering
                    {\unitlength=1cm
                    \begin{picture}(12,7.5)(0,-0.5)
                    \put(1,0){\includegraphics[width=0.7\textwidth,clip]{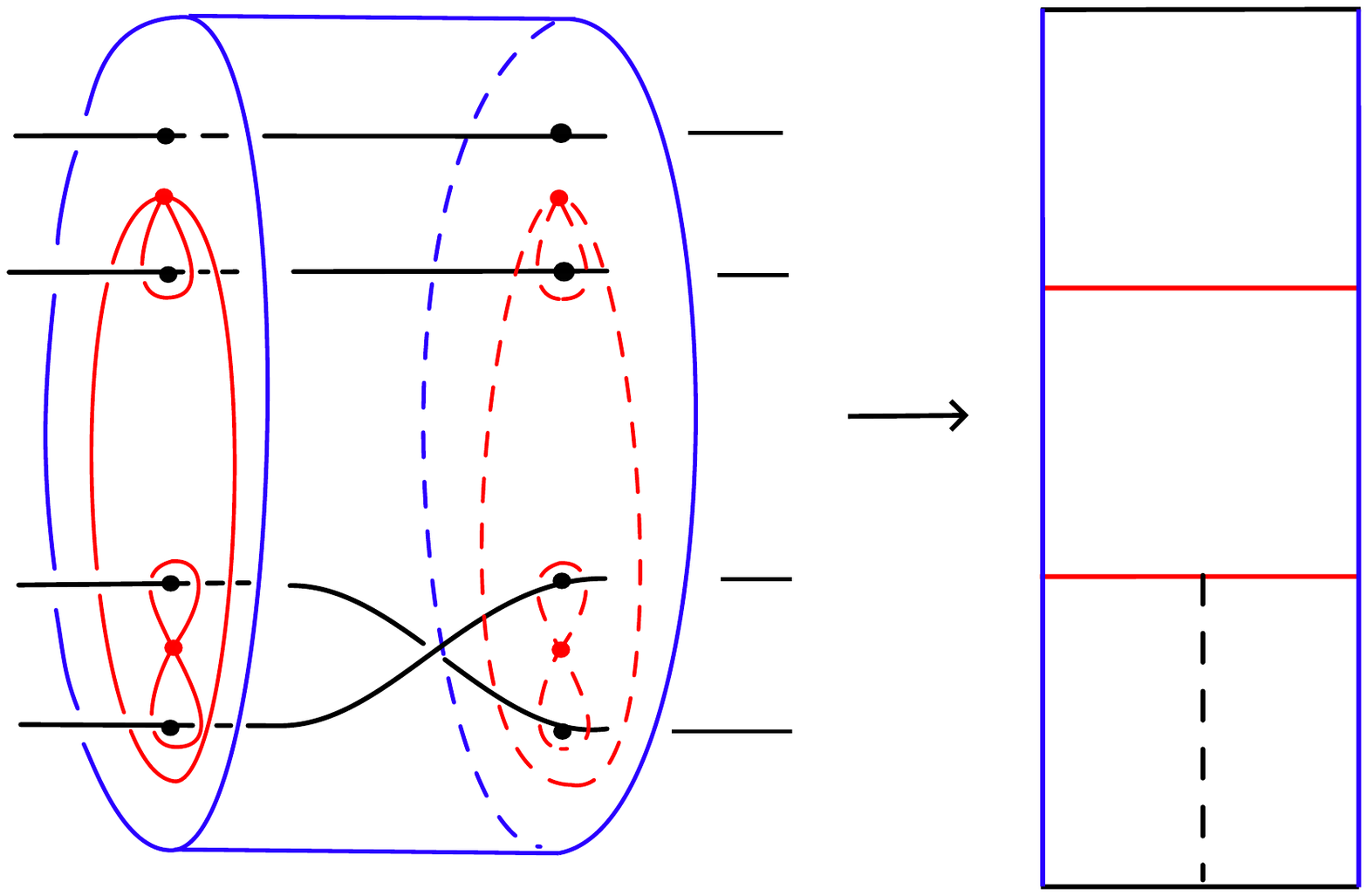}}
                    \put(1,6.3){\underline{\textbf{Type $3$}}}
                    \put(6.65,3.45){$\Phi_{k}$}
                    \put(1.3,5.5){$N_{k}$}
                    \put(7,5.5){$T_{k}$}

                    \put(1,5){$1$}
                    \put(1,4){$2$}
                    \put(1,2.2){$3$}
                    \put(1,1.3){$4$}

                    \put(2,0.1){$F_{k}$}
                    \put(4.8,0.1){$F_{k+1}$}
                    \put(7.5,0){$\gamma_{k}$}
                    \put(9.5,0){$\gamma_{k+1}$}
                    \end{picture}}
                    \caption{ A smooth map to the region of Type $3$. }
                    \label{type32}
                \end{figure}

For $N_k$ and the rectangular region $T_k$ of Type 3, in a similar way as before, we construct a smooth map $\Phi_{k} : N_{k} \to T_{k}$ as natural extensions of $\psi_{k}$ and $\psi_{k+1}$ by identifying $F_k$ and $F_{k+1}$ with $N_{k-1} \cap N_{k}$ and $N_{k} \cap N_{k+1}$, respectively. 
See Fig.~\ref{type32}. 
Also, this map $\Phi_{k}$ is the same as that given in the previous section shown in Fig.~\ref{type3}, but we include Fig.~\ref{type32} for consistency with the deformation described above. 
In this case, we note that there are no fibers containing two indefinite fold points. 

Now, we connect $N_{0}, N_{1}, \cdots, N_{n-1}$ and $N_{n}$ by identifying the boundaries of them to obtain $S^{3}$. 
Also, we connect $T_{0}, T_{1}, \cdots, T_{n-1}$ and $T_{n}$ by identifying $\gamma_{1}, \cdots, \gamma_{n}$ to make $E$. 
Then, a smooth map $f_{3} : S^{3} \to \mathbb{R}^{2}$ satisfying $f_{3}|_{N_{k}} = \Phi_{k}$ for $k \in \{ 0,1,2,\cdots,m \}$ is constructed. 
By construction, this map $f_3$ is a stable map from $S^3$ into $\mathbb{R}^{2}$ with desired properties. 
In particular, $f_3$ has no singular fibers of type $\mathrm{I\hspace{-1.2pt}I^{2}}$ and the total number of singular fibers of type $\mathrm{I\hspace{-1.2pt}I^{3}}$ is $\frac{1}{2} \sum_{i=1}^{m} | b_{i} |$. 
This completes a proof of Theorem~\ref{main2}. 

        \begin{example}
        A stable map $f : S^{3} \to \mathbb{R}^2$ with the Whitehead link as $S_{0}(f)$ obtained by the construction given in this section is presented in Fig.~\ref{$49$}. 
        This $f$ has just one fiber of type $\mathrm{I\hspace{-1.2pt}I^{3}}$. 
        Some of the other fibers containing an indefinite fold point are also shown in the figure. 

        \begin{figure}[H]
                \setlength\unitlength{1truecm}
                \begin{picture}(12.5,12.5)(0,0)
                    \put(1,0){\includegraphics[width=0.75\textwidth,clip]{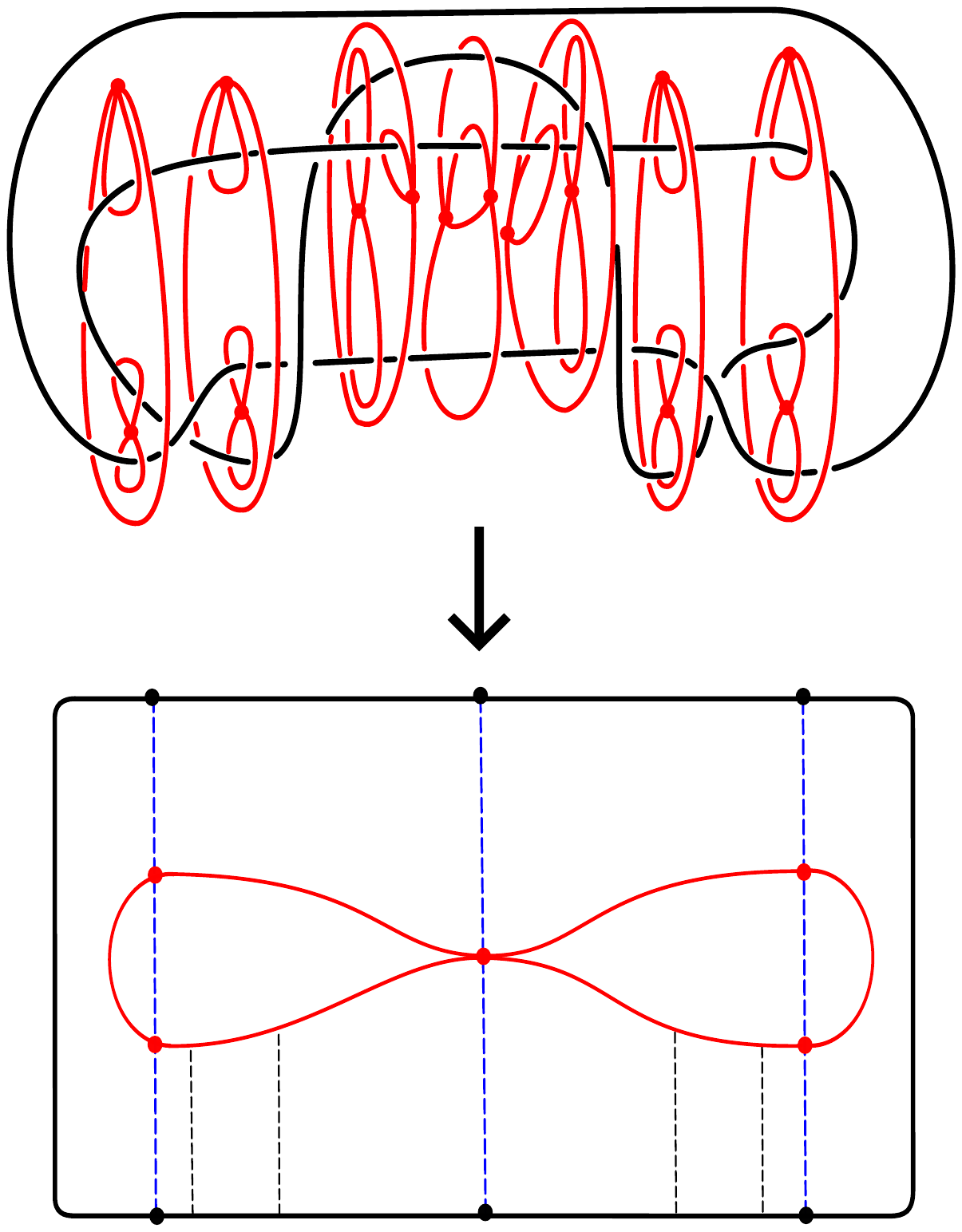}}
                \end{picture}
                \caption{A stable map $f$ from $S^{3}$ with the Whitehead link as $S_{0}(f)$.}
                \label{$49$}
        \end{figure}

        \end{example}

 \section{Stable map complexities}

In \cite[Section 6]{Ishikawa-Koda}, the relationship between stable map complexities and hyperbolic volumes for 3-manifolds is discussed.
In this section, based on that, we give an estimate on stable map complexities of two-bridge link exteriors by using the constructions given in the previous sections.

In \cite{Ishikawa-Koda}, the following relationship between the hyperbolic volumes and the stable map complexities is obtained.
Let $M$ be a compact 3-manifold $M$ with (possibly empty) boundary consisting of tori.
Suppose that the interior $\mathrm{int}M$ of $M$ admits a complete hyperbolic metric of finite volume.
Then,
\begin{equation}\label{eq41}
\mathrm{vol} (M) \le 2 V_{oct}\,\mathrm{smc}(M)    
\end{equation}
holds, where $\mathrm{vol}(M)$ denotes the hyperbolic volume of $\mathrm{int}M$, $\mathrm{smc}(M)$ the \textit{stable map complexity} of $M$, that is, the minimal number of the weighted sums
$ \left| \mathrm{I\hspace{-1.2pt}I^{2}}(f) \right| + 2 \left| \mathrm{I\hspace{-1.2pt}I^{3}}(f) \right| $
for stable maps $f : M \to \mathbb{R}^2$, and $V_{oct}$ the volume of the hyperbolic ideal regular octahedron, which is approximately $3.6638$.

Let $L$ be a two-bridge link in $S^3$ which admits a Conway form $C(a_{1},b_{1},\cdots,a_{m},b_{m},a_{m+1} )$ with all the $b_{i}$'s are even.
Then, by Theorem~\ref{main1}, we have a stable map of $S^3$ to $\mathbb{R}^2$ with $|\mathrm{I\hspace{-1.2pt}I^{2}}(f_{2})| = 2m$ and $\mathrm{I\hspace{-1.2pt}I^{3}}(f_{2}) = \emptyset$.
From this map, by removing an open tubular neighborhood $N(L)$ of $L$ from $S^3$ suitably, a stable map $f : E(L) \to \mathbb{R}^2$ is obtained for the exterior $E(L) = S^3 \setminus N(L)$.
This map $f$ satisfies that $ \left| \mathrm{I\hspace{-1.2pt}I^{2}}(f) \right| = 2m$ and $\mathrm{I\hspace{-1.2pt}I^{3}}(f) = \emptyset$, and so,
\[\left| \mathrm{I\hspace{-1.2pt}I^{2}}(f) \right| + 2 \left| \mathrm{I\hspace{-1.2pt}I^{3}}(f) \right| = 2m . \]
This gives
\begin{equation}\label{eq42}
\mathrm{smc}( E(L) ) \le 2m
\end{equation}
for the two-bridge link $L$.

    On the other hand, if the Conway form $C(a_{1},b_{1},\cdots,a_{m},b_{m},a_{m+1} )$ is a reduced alternating diagram, equivalently, all the $a_i$'s and $b_j$'s are at least 2, then the following holds by \cite[Theorem B.3]{FuterGueritaud06}.
    Let $D=C(a_{1},b_{1},\cdots,a_{m},b_{m},a_{m+1} )$ be a reduced alternating diagram of a hyperbolic two-bridge link $L$.
    Then
    \begin{equation}\label{eq43}
    \mathrm{vol} ( S^3 \setminus L) < 2(\mathrm{tw} (D) - 1) v_{oct}
    \end{equation}
    and actually $\mathrm{tw}(D) = 2m+1$.
    Here, the \textit{twist number} $\mathrm{tw}(D)$ denotes the number of equivalence classes of crossings, where two crossings are considered equivalent if there is a loop in the projection plane intersecting $D$ transversely precisely in the two crossings.
    
    Moreover, as remarked in \cite[Appendix B]{FuterGueritaud06}, if $a_{i}$ and $b_{j}$ are sufficiently large, then the volume $\mathrm{vol} ( S^3 \setminus L)$ is close enough to the number $2(\mathrm{tw} (D) - 1) v_{oct} = 4m v_{oct}$.
    
By combining Equations~\eqref{eq41}, \eqref{eq42}, and \eqref{eq43}, if $a_{i}$ and $b_{j}$ are sufficiently large, we have 
\begin{equation*}
4m v_{oct} - \varepsilon < \mathrm{vol} (S^{3} \setminus L) < 2v_{oct} \mathrm{smc}(E(L)) \leq 4m v_{oct}
\end{equation*}
for arbitrary small $\varepsilon > 0$.

This implies that $\mathrm{smc}(E(L)) = 2m$ holds for 
a two-bridge link $L$ which has the Conway form $C(a_{1},b_{1},\cdots,a_{m},b_{m},a_{m+1})$ with $a_{i}$, $b_{j}$ sufficiently large even. 
This completes a proof of Corollary~\ref{cor}. 



\begin{remark}
The following were suggested by the anonymous referee: 
There are several complexities for links. 
The most famous complexity is the minimum number of ideal tetrahedrons to realize the complement of a link. 
This complexity is studied in many papers. 
For example, recently in \cite{morgan2024complexity2bridgelinkcomplements}, 
the complexity of two-bridge links and their relationship to the Sakuma-Weeks triangulation are studied. 
Studying the complexity of two-bridge links is an interesting topic in recent research. 
The results in this paper can be seen as a variation of such studies of complexities. 
\end{remark}

\section*{Acknowledgements}
The authors would like to thank Masaharu Ishikawa for giving us an idea about the construction of the stable map $f_{2}$ in Theorem \ref{main1}. They also thank Osamu Saeki, Takahiro Yamamoto, Yuya Koda, and Ryoga Furutani for their helpful advice. The first author is partially supported by JSPS KAKENHI Grant Number JP22K03301.

\bibliographystyle{ws-jktr}
\bibliography{name}

\end{document}